\documentclass[11pt]{article}
\usepackage{setspace}
 
\usepackage{amsmath} 
\usepackage{amssymb}
\usepackage{latexsym}
\usepackage{xcolor}
\usepackage{fancyhdr}
\allowdisplaybreaks 

 \usepackage{comment}

\def\XXint#1#2#3{{\setbox0=\hbox{$#1{#2#3}{\int}$} 
\vcenter{\hbox{$#2#3$}}\kern-.5\wd0}}   

 \numberwithin{equation}{section}
\newtheorem{theorem}[equation]{Theorem}
\newtheorem{proposition}[equation]{Proposition}
\newtheorem{definition}[equation]{Definition}
\newtheorem{remark}[equation]{Remark}
\newtheorem{lemma}[equation]{Lemma}
\newtheorem{corollary}[equation]{Corollary}
\newtheorem{example}[equation]{Example}
\title{H\"{o}lder regularity of the solutions of 
Fredholm integral equations  on upper Ahlfors regular sets} 
 
 \author{  
Massimo Lanza de Cristoforis\thanks{
  Dipartimento di Matematica `Tullio Levi-Civita', 
Universit\`a degli Studi di Padova, 
Via Trieste 63, Padova 35121, 
Italy. E-mail: mldc@math.unipd.it}  
and Manuel Norman \thanks{Dipartimento di Matematica `Tullio Levi-Civita', 
Universit\`a degli Studi di Padova, 
Via Trieste 63, Padova 35121, 
Italy. E-mail: manuel.norman@studenti.unipd.it}
}

\date{\ }

\begin{document}
 
 \maketitle

\noindent
{\bf Abstract:}  We extend to the context of metric measured spaces, with a measure that satisfies upper Ahlfors growth conditions  the validity of    (generalized) H\"{o}lder continuity results for the solution of a Fredholm integral equation of the second kind. Here we note that upper Ahlfors growth conditions  include also cases of nondoubling measures.
 \vspace{\baselineskip}

\par
 

\section{Introduction} This paper concerns the treatment of Fredholm integral equations of the second kind on subsets $Y$ of the Euclidean space ${\mathbb{R}}^n$ with a measure.
A typical example is the case in which $Y$ equals the boundary of a bounded open subset of ${\mathbb{R}}^n$ (so for example the so-called case of the boundary integral equations) or the case in which $Y$  is a manifold of codimension greater than or equal to $1$ that is embedded into ${\mathbb{R}}^n$. The goal here is to provide general statements that could be applied in a whole variety of cases and the plan is to   present a unified approach to some results by assuming that   $Y$ is a subset  of a metric space $(M,d)$ and that $Y$ is equipped with a measure $\nu$ that satisfies upper Ahlfors growth conditions that include non-doubling measures that we introduce below.\par
 
Since the analysis of integral equations involves the analysis of integral operators, we mention that Edmunds, Kokilashvili and Meskhi \cite{EdKoMe02}, Garcia-Cuerva and Gatto \cite{GaGa04}, \cite{GaGa05}, Gatto \cite{Gat06}, \cite{Gat09}, Kokilashvili, Meskhi, Rafeiro and Samko \cite{KoMeRaSa16}, Mattila, Melnikov and Verdera \cite{MaMeVe96}, Mitrea, Mitrea, Mitrea \cite{MitMitMit22_6}, Nazarov, Treil, Volberg \cite{NaTrVo97}, Verdera \cite{Ve13}, \cite{Ve23}  and Tolsa \cite{To98}, \cite{To01a}, \cite{To01b} have worked on the theory of integral operators acting in metric spaces with a measure.\par

Here we consider a subset $Y$ of $M$ and an   auxiliary subset $X$  of $M$. Moreover, we formulate the following assumption.
\begin{eqnarray} \nonumber
&&\text{Let}\ {\mathcal{N}} \ \text{ be a $\sigma$-algebra of parts of}\  Y  \,, {\mathcal{B}}_Y\subseteq  {\mathcal{N}}\,.
\\   \label{eq:nu}
&&\text{Let}\ \nu\  \text{be   measure on}\  {\mathcal{N}} \,.
\\  \nonumber
&&\text{Let}\ \nu(B(x,r)\cap Y)<+\infty\qquad\forall (x,r)\in X\times ]0,+\infty[\,.
\end{eqnarray}
Here  ${\mathcal{B}}_Y$ denotes the $\sigma$-algebra of the Borel subsets of $Y$ and 
\begin{equation}\label{eq:balls}
B(\xi,r)\equiv \left\{\eta\in M:\, d(\xi,\eta)<r\right\}\,,\quad
B(\xi,r]\equiv \left\{\eta\in M:\, d(\xi,\eta)\leq r\right\}
\,,
\end{equation}
for all $(\xi,r)\in M\times ]0,+\infty[$. We assume that  $\upsilon_Y\in ]0,+\infty[$ and we consider two types of assumptions on $\nu$. 
The first assumption is  that
$Y$ is   upper $\upsilon_Y$-Ahlfors regular with respect to $X$, \textit{i.e.}, that
\begin{eqnarray} \nonumber
&&\text{there\ exist}\ r_{X,Y,\upsilon_Y}\in]0,+\infty]\,,\ c_{X,Y,\upsilon_Y}\in]0,+\infty[\ \text{such\ that}
\\ \nonumber
&&\nu( B(x,r)\cap Y )\leq c_{X,Y,\upsilon_Y} r^{\upsilon_Y}  
\\  \label{defn:uareg1}
&&\text{for\ all}\ x\in X\ \text{and}\  r \in]0,r_{X,Y,\upsilon_Y}[
 \,. 
\end{eqnarray}
In case $X=Y$, we just say that $Y$ is   upper $\upsilon_Y$-Ahlfors regular  and this is the assumption that has been considered by  
Garc\'{\i}a-Cuerva and Gatto \cite{GaGa04}, \cite{GaGa05}, Gatto \cite{Gat06}, \cite{Gat09}   in  the context of H\"{o}lder spaces. See also  Edmunds,  Kokilashvili and   Meskhi~\cite[Chap.~6]{EdKoMe02} in the frame of Lebsgue spaces. 

An interesting feature of condition (\ref{defn:uareg1}) is that it does not imply any estimate  of $\nu( B(x,r)\cap Y )$ from below in terms of $r^{\upsilon_Y}$ as in the so called lower $\upsilon_Y$-Ahlfors regularity condition  that together with (\ref{defn:uareg1}) would imply the validity
of the so-called $\upsilon_Y$-Ahlfors regularity condition (as in David and Semmes \cite{DaSe93}) and accordingly the validity of the so-called doubling  condition for the measure $\nu$.
However, as pointed out in the papers of Verdera \cite{Ve13}, \cite[p.~21]{Ve23}  in connection with the integral operator of the Cauchy kernel and later in  \cite[Example 6.1]{La24}, both the doubling and the Ahlfors conditions may be too stringent. 

As in \cite{La22a}, we also  consider a stronger condition than (\ref{defn:uareg1}) that still does not involve estimates from below for $\nu$  in which we replace the ball $B(x,r)\cap Y$ in $Y$ with an annular domain 
$(B(x,r_2)\setminus B(x,r_1))\cap Y $ with $0\leq r_1< r_2$ in $Y$ and that  says that  one can estimate the $\nu$-measure of an annular domain 
$(B(x,r_2)\setminus B(x,r_1))\cap Y$ in $Y$ from above  in terms of the measure of an annular domain of radii $r_1$ and $r_2$ in a Euclidean space of dimension $\upsilon_Y$ (at least for integer values of $\upsilon_Y$).  Namely, 
  we assume that $Y$ is strongly upper $\upsilon_Y$-Ahlfors regular with respect to $X$, \textit{i.e.}, that
\begin{eqnarray*} \nonumber
&&\text{there\ exist}\ r_{X,Y,\upsilon_Y}\in]0,+\infty]\,,\ c_{X,Y,\upsilon_Y}\in]0,+\infty[\ \text{such\ that}
\\ \nonumber
&&\nu( (B(x,r_2)\setminus B(x,r_1))\cap Y )\leq c_{X,Y,\upsilon_Y}(r_2^{\upsilon_Y}-r_1^{\upsilon_Y})
\\ \nonumber 
&&\text{for\ all}\ x\in X\ \text{and}\ r_1,r_2\in[0,r_{X,Y,\upsilon_Y}[
\ \text{with}\ r_1<r_2\,,
\end{eqnarray*}
where we understand that $B(x,0)\equiv\emptyset$  (in case $X=Y$, we just say that $Y$ is strongly upper $\upsilon_Y$-Ahlfors regular). So, for example,	
 if $Y$ is the boundary of a bounded open Lipschitz  subset of $M={\mathbb{R}}^n$, then $Y$ is upper $(n-1)$-Ahlfors regular with respect to ${\mathbb{R}}^n$ (cf.~
 \cite[Prop.~6.5]{La24}) 
and if $Y$ is the boundary of an open bounded subset of $M={\mathbb{R}}^n$ of class $C^1$, then $Y$ is strongly upper $(n-1)$-Ahlfors regular  with respect to $Y$ (cf.~
\cite[Prop.~6.6, Rem.~6.2]{La24}).

In the present paper, we assume that if $Y$  (strongly) upper $\nu$-Ahlfors regular with respect to
$Y$, and we consider the problem of whether
a classical   $L^2_\nu(Y)$-solution  of a Fredholm integral equation of the second type with either continuous or H\"{o}lder continuous data actually carries the same regularity of the data. More precisely, we set
\[
{\mathbb{D}}_{Y\times Y}\equiv \{(y_1,y_2)\in Y\times Y:\,y_1=y_2\}\,,
\]
we consider a continuous kernel $K$ from $(Y\times Y)\setminus {\mathbb{D}}_{Y\times Y}$ to ${\mathbb{C}}$, we assign a bounded and continuous (or H\"{o}lder continuous) datum $g$ on $Y$, 
we assume that the Fredholm integral equation
\begin{equation}\label{eq:frsek}
\mu- A[K,\mu]=g \qquad\text{on}\ Y\,,
\end{equation}
where
\begin{equation}
\label{prop:k0a0}
A[K,\mu](x)\equiv\int_YK(x,y)\mu(y)\,d\nu(y)\qquad\forall x\in X
\end{equation}
has a solution $\mu\in L_\nu^2(Y)$ and we ask whether $\mu$ would also be continuous (or H\"{o}lder continuous). Questions of this type are well known in the classical theory in case
$Y$ is either a measurable subset of ${\mathbb{R}}^n$ with the ordinary Lebesgue measure
or the boundary of a bounded open Lipschitz subset of  ${\mathbb{R}}^n$ with the ordinary surface measure, thus with a measure that satisfies the (double sided) Ahlfors regularity condition.

In  the present paper we plan to extend such results to the context of metric measured spaces, with measures that satisfy only upper  Ahlfors regularity conditions (see Theorems \ref{thm:cizere}, \ref{thm:cizereh}, \ref{thm:cizereht}).

Here we note that the  Fredholm theory has important applications to differential equations, since many boundary value problems on open domains can be reduced to the study of certain Fredholm integral equations (see \textit{e.g.},  Folland~\cite{Fo95}).

Many results have been shown to hold for integral equations on the boundary of   a bounded open Lipschitz subset of  ${\mathbb{R}}^n$, and thus in case the measure   satisfies a (double sided) Ahlfors regularity condition.   However, there are again cases where the Ahlfors condition should be replaced by the upper Ahlfors condition, as  with cusp domains (cf., \textit{e.g.}, \cite[Ex.~6.1]{La24}).

The paper is organized as follows. In Section \ref{sec:tecprel}, we introduce our notation and some basic definitions.   
In section \ref{sec:kecla},  we introduce some properties of potential type kernels and some technical results in metric measured spaces. 
In section \ref{sec:cmpoke}, we extend some classical results on composite kernels to our metric measured space context.

In Section \ref{sec:corefr}, we prove that if $g$ is continuous and bounded, then a solution $\mu\in L_\nu^2(Y)$ of the Fredholm equation (\ref{eq:frsek}) is also continuous and bounded. In Section \ref{sec:horefr}, we prove that if $g$ is H\"{o}lder continuous and bounded, then a solution $\mu\in L_\nu^2(Y)$ of the Fredholm equation (\ref{eq:frsek}) satisfies a condition of generalized  H\"{o}lder  continuity.
  
\section{Basic notation}\label{sec:tecprel}
Let $X$ be a set. Then we set
\[
B(X)\equiv\left\{
f\in {\mathbb{C}}^X:\,f\ \text{is\ bounded}
\right\}
\,,\quad
\|f\|_{B(X)}\equiv\sup_X|f|\qquad\forall f\in B(X)\,.
\]
 Then $C^0(M)$ denotes the set of continuous functions from $M$ to ${\mathbb{C}}$ and we introduce the subspace
$
C^0_b(M)\equiv C^0(M)\cap B(M)
$
of $B(M)$.  Let $\omega$ be a function from $[0,+\infty[$ to itself such that
\begin{eqnarray}
\nonumber
&&\qquad\qquad\omega(0)=0,\qquad \omega(r)>0\qquad\forall r\in]0,+\infty[\,,
\\
\label{om}
&&\qquad\qquad\omega\ {\text{is\   increasing,}}\ \lim_{r\to 0^{+}}\omega(r)=0\,,
\\
\nonumber
&&\qquad\qquad{\text{and}}\ \sup_{(a,t)\in[1,+\infty[\times]0,+\infty[}
\frac{\omega(at)}{a\omega(t)}<+\infty\,.
\end{eqnarray}
If $f$ is a function from a subset ${\mathbb{D}}$ of a metric space $(M,d)$   to ${\mathbb{C}}$,  then we denote by   $|f:{\mathbb{D}}|_{\omega (\cdot)}$  the $\omega(\cdot)$-H\"older constant  of $f$, which is delivered by the formula   
\[
|f:{\mathbb{D}}|_{\omega (\cdot)
}
\equiv
\sup\left\{
\frac{|f( x )-f( y)|}{\omega(d( x, y))
}: x, y\in {\mathbb{D}} ,  x\neq
 y\right\}\,.
\]        
If $|f:{\mathbb{D}}|_{\omega(\cdot)}<+	\infty$, we say that $f$ is $\omega(\cdot)$-H\"{o}lder continuous. Sometimes, we simply write $|f|_{\omega(\cdot)}$  
instead of $|f:{\mathbb{D}}|_{\omega(\cdot)}$. The
subset of $C^{0}({\mathbb{D}} ) $  whose
functions  are
$\omega(\cdot)$-H\"{o}lder continuous    is denoted  by  $C^{0,\omega(\cdot)} ({\mathbb{D}})$
and $|f:{\mathbb{D}}|_{\omega(\cdot)}$ is a semi-norm on $C^{0,\omega(\cdot)} ({\mathbb{D}})$.  
Then we consider the space  $C^{0,\omega(\cdot)}_{b}({\mathbb{D}} ) \equiv C^{0,\omega(\cdot)} ({\mathbb{D}} )\cap B({\mathbb{D}} ) $ with the norm \[
\|f\|_{ C^{0,\omega(\cdot)}_{b}({\mathbb{D}} ) }\equiv \sup_{x\in {\mathbb{D}} }|f(x)|+|f|_{\omega(\cdot)}\qquad\forall f\in C^{0,\omega(\cdot)}_{b}({\mathbb{D}} )\,.
\] 
 \begin{remark}
\label{rem:om4}
Let $\omega$ be as in (\ref{om}). 
Let ${\mathbb{D}}$ be a   subset of ${\mathbb{R}}^{n}$. Let $f$ be a bounded function from $ {\mathbb{D}}$ to ${\mathbb{C}}$, $a\in]0,+\infty[$.  Then,
\[
\label{rem:om5}
\sup_{x,y\in {\mathbb{D}},\ |x-y|\geq a}\frac{|f(x)-f(y)|}{\omega(|x-y|)}
\leq \frac{2}{\omega(a)} \sup_{{\mathbb{D}}}|f|\,.
\]
\end{remark}
In the case in which $\omega(\cdot)$ is the function 
$r^{\alpha}$ for some fixed $\alpha\in]0,1]$, a so-called H\"{o}lder exponent, we simply write $|\cdot:{\mathbb{D}}|_{\alpha}$ instead of
$|\cdot:{\mathbb{D}}|_{r^{\alpha}}$, $C^{0,\alpha} ({\mathbb{D}})$ instead of $C^{0,r^{\alpha}} ({\mathbb{D}})$, $C^{0,\alpha}_{b}({\mathbb{D}})$ instead of $C^{0,r^{\alpha}}_{b} ({\mathbb{D}})$, and we say that $f$ is $\alpha$-H\"{o}lder continuous provided that 
$|f:{\mathbb{D}}|_{\alpha}<+\infty$.   For each $\theta\in]0,1]$, we define the function $\omega_{\theta}(\cdot)$ from $[0,+\infty[$ to itself by setting
\begin{equation}
\label{omth}
\omega_{\theta}(r)\equiv
\left\{
\begin{array}{ll}
0 &r=0\,,
\\
r^{\theta}|\ln r | &r\in]0,r_{\theta}]\,,
\\
r_{\theta}^{\theta}|\ln r_{\theta} | & r\in ]r_{\theta},+\infty[\,,
\end{array}
\right.
\end{equation}
where
$
r_{\theta}\equiv e^{-1/\theta}
$ for all $\theta\in ]0,1]$. Obviously, $\omega_{\theta}(\cdot) $ is concave and satisfies   condition (\ref{om}).

We retain the standard notation for the Lebesgue space $L^{p}_\nu(Y)$  of $p$-summable (equivalence classes of) measurable functions for $p\in[1,+\infty]$. In particular, we set 
\[
\|f\|_{L^{p}_\nu(Y)}\equiv \left(\int_{Y}|f|^p\,d\nu\right)^{1/p}\qquad\forall f\in L^{p}_\nu(Y) 
\]
if $p\in [1,+\infty[$ and
\[
\|f\|_{L^{\infty}_\nu(Y)}\equiv {\mathrm{ess\ sup}}_Y|f|\qquad\forall f\in L^{\infty}_\nu(Y) \,.
\]
For basic inequalities such as the H\"{o}lder inequality for Lebesgue spaces, we refer to textbooks such as Folland~\cite[Chap.~6]{Fo99}.

 \section{Special classes of potential type kernels in metric spaces}\label{sec:kecla}
 If $X$ and $Y$ are sets, then we denote by ${\mathbb{D}}_{X\times Y}$ the diagonal  of $X\times Y$, i.e., we set
 \begin{equation}\label{diagonal}
{\mathbb{D}}_{X\times Y}\equiv\left\{
(x,y)\in X\times Y:\,x=y
\right\}\,. 
\end{equation}
An off-diagonal		 function  in $X\times Y$ is a function from $(X\times Y)\setminus {\mathbb{D}}_{X\times Y}$ to ${\mathbb{C}}$. We now wish to consider a specific class of off-diagonal kernels in a metric space  $(M,d)$.  
\begin{definition}
 Let $X$ and $Y$ be subsets of a metric space  $(M,d)$. Let $s\in {\mathbb{R}}$. We denote by ${\mathcal{K}}_{s,X\times Y}$, the set of continuous functions $K$ 
 from $ (X\times Y)\setminus {\mathbb{D}}_{ X\times Y }$ to ${\mathbb{C}}$ such that
\[
 \|K\|_{ {\mathcal{K}}_{s,X\times Y} }\equiv \sup_{(x,y)\in  (X\times Y)\setminus {\mathbb{D}}_{ X\times Y }  }
 |K(x,y)|\,d(x,y)^s<+\infty\,.
\]
The elements of $ {\mathcal{K}}_{s,X\times Y}$ are said to be kernels of potential type $s$ in $X\times Y$. 
\end{definition}
We plan to consider `potential type' kernels as in the following definition (see also paper \cite{DoLa17} with Dondi,  where such classes have been introduced in a form that generalizes those of Giraud \cite{Gi34}, Gegelia \cite{Ge67}, 
 Kupradze, Gegelia, Basheleishvili and 
 Burchuladze \cite[Chapter~IV]{KuGeBaBu79}).	 
\begin{definition}\label{defn:ksss}
 Let $X$ and $Y$ be subsets of a metric space  $(M,d)$. Let $s_1$, $s_2$, $s_3\in {\mathbb{R}}$. We denote by ${\mathcal{K}}_{s_1, s_2, s_3} (X\times Y)$ the set of continuous functions $K$ from $(X\times Y)\setminus {\mathbb{D}}_{X\times Y}$ to ${\mathbb{C}}$ such that
 \begin{eqnarray*}
\lefteqn{
\|K\|_{  {\mathcal{K}}_{ s_1, s_2, s_3  }(X\times Y)  }
\equiv
\sup\biggl\{\biggr.
d(x,y)^{ s_{1} }|K(x,y)|:\,(x,y)\in X\times Y, x\neq y
\biggl.\biggr\}
}
\\ \nonumber
&&\qquad\qquad\qquad
+\sup\biggl\{\biggr.
\frac{d(x',y)^{s_{2}}}{d(x',x'')^{s_{3}}}
|  K(x',y)- K(x'',y)  |:\,
\\ \nonumber
&&\qquad\qquad\qquad 
x',x''\in X, x'\neq x'', y\in Y\setminus B(x',2d(x',x''))
\biggl.\biggr\}<+\infty\,.
\end{eqnarray*}
\end{definition}
 One can easily verify that $({\mathcal{K}}_{ s_{1},s_{2},s_{3}   }(X\times Y),\|\cdot\|_{  {\mathcal{K}}_{s_{1},s_{2},s_{3}   }(X\times Y)  })$ is a normed space.  By our definition, if $s_1$, $s_2$, $s_3\in {\mathbb{R}}$, we have
\[
{\mathcal{K}}_{s_{1},s_{2},s_{3}   }(X\times Y) \subseteq {\mathcal{K}}_{s_{1}, X\times Y} 
\]
and
\[
\|K\|_{{\mathcal{K}}_{s_{1}, X\times Y} }\leq \|K\|_{ {\mathcal{K}}_{s_{1},s_{2},s_{3}   }(X\times Y) }
\qquad\forall K\in {\mathcal{K}}_{s_{1},s_{2},s_{3}   }(X\times Y) \,.
\]
We note that if we choose $s_2=s_1+s_3$ we have a so-called class of standard kernels. 
\begin{remark}\label{rem:bddimfi}
 Let  $X$ and  $Y$ be subsets
  of a metric space  $(M,d)$.    Let $\nu$ be as in (\ref{eq:nu}) .   If $Y\neq \emptyset$, then we further require that 
  $X\neq\emptyset$.	   If $Y$ is bounded, then $\nu(Y)<+\infty$. Indeed, if $Y=\emptyset$, then $\nu(Y)=0$. If instead $Y\neq\emptyset$, then there exist at least an element 
  $\tilde{x}\in X$ and an element	 	$\tilde{y}\in Y$ and the assumption that  $Y$ is bounded 
  together with the triangular inequality imply	 
  that there exists $r\in]0,+\infty[$ such that $Y\subseteq B(\tilde{x},r)$ and accordingly condition (\ref{eq:nu}) implies that $\nu(Y)$ is finite. 

\end{remark}

Then one can prove the following  basic inequalities  for the integral on an upper Ahlfors    regular set $Y$ and on the intersection of $Y$ with balls with center at a point $x$ of $X$ 
of the powers of $d(x,y)^{-1}$ with exponent  $s\in]-\infty,\upsilon_Y[$, that are variants of those proved by Gatto \cite[p.~104]{Gat09} in case $X=Y$ (for a proof see \cite[Lem.~3.2, 3.4]{La22a}).
\begin{lemma}\label{lem:cominur} Let $X$ and $Y$ be subsets of a metric space  $(M,d)$.
 Let $\nu$ be as in (\ref{eq:nu}).  Let $\upsilon_Y\in]0,+\infty[$.  Let  $Y$ be  upper $\upsilon_Y$-Ahlfors regular with respect to $X$. Then the following statements hold.
\begin{enumerate}
\item[(i)] $\nu(\{x\})=0$ for all $x\in X\cap Y$.
\item[(ii)] Let $\nu(Y)<+\infty$. If $s\in ]0,\upsilon_Y[$, then
\[
c'_{s,X,Y}\equiv\sup_{x\in X}\int_Y\frac{d\nu(y)}{d(x,y)^s}
\leq \nu(Y)a^{-s}+c_{X,Y,\upsilon_Y}\frac{\upsilon_Y}{\upsilon_Y-s}a^{\upsilon_Y-s}
\]
for all $a\in]0,r_{X,Y,\upsilon_Y}[$. If $s=0$, then 
\[
c'_{0,X,Y}\equiv\sup_{x\in X}\int_Y\frac{d\nu(y)}{d(x,y)^0}
= \nu(Y) \,.
\]
\item[(iii)]  Let  $\nu(Y)<+\infty$ whenever $r_{X,Y,\upsilon_Y}<+\infty$. If  $s\in]-\infty,\upsilon_Y[$, then
\[
c''_{s,X,Y}\equiv\sup_{(x,t)\in X\times]0,+\infty[}
t^{s-\upsilon_Y}\int_{B(x,t)\cap Y}\frac{d\nu(y)}{d(x,y)^s}
<+\infty
\,.
\]
\end{enumerate}
\end{lemma}

Then one can prove the following  basic inequalities  for the integral  of the powers of $d(x,y)^{-1}$ with exponent  $s\in[\upsilon_Y,+\infty[$
on the complement in $Y$ of balls with center at a point $x$ of $X$,   that are variants of those proved by Gatto \cite[p.~104]{Gat09} in case $X=Y$ (for a proof see \cite[Lem.~3.6]{La22a}).
  \begin{lemma}\label{lem:cominu} Let 
$
(M,d)$ be a metric space. Let $X$, $Y\subseteq M$. 
  Let $\nu$ be as in (\ref{eq:nu}). Let  $\nu(Y)<+\infty$. Let $\upsilon_Y\in]0,+\infty[$. Then the following statements hold.
\begin{enumerate}
\item[(i)] Let  $Y$ be   upper $\upsilon_Y$-Ahlfors regular with respect to $X$.  If $s\in ]\upsilon_Y,+\infty[$, then 
\[
c'''_{s,X,Y}\equiv\sup_{(x,t)\in X\times]0,+\infty[}
t^{s-\upsilon_Y}\int_{Y\setminus B(x,t) }\frac{d\nu(y)}{d(x,y)^s}<+\infty\,.
\]
\item[(ii)] Let  $Y$ be strongly upper $\upsilon_Y$-Ahlfors regular with respect to $X$.   Then
\[
c^{iv}_{X,Y}\equiv\sup_{(x,t)\in X\times]0,1/e[}
\vert\ln t\vert^{-1}\int_{Y\setminus B(x,t) }\frac{d\nu(y)}{d(x,y)^{\upsilon_Y}}<+\infty\,.
\]
\end{enumerate}
\end{lemma}
We also note that the following statement holds true (cf.~\cite[Lem.~4.2]{La22b}).
\begin{lemma}\label{lem:cominuu} Let $X$ and $Y$ be subsets of a metric space  $(M,d)$. Let $\nu$ be as in (\ref{eq:nu}). Let $\nu(Y)<+\infty$. Let $\upsilon_Y\in]0,+\infty[$,  $s\in [0,\upsilon_Y[$.    
 Let  $Y$ be  upper $\upsilon_Y$-Ahlfors regular with respect to $X$. Then for each  $\epsilon\in]0,+\infty[$  there exists $\delta\in ]0,+\infty[$ such that 
\begin{equation} \label{lem:cominuu1}
\sup_{x\in X}\int_{F }
 d(x,y)^{-s} \,d\nu(y) \leq\epsilon
\qquad \text{if}\ F\in {\mathcal{N}},\ \nu(F)\leq \delta\,.
\end{equation}
\end{lemma}
Finally,  we introduce the following  known 	elementary lemma, which we exploit later and which can be proved by the triangular inequality.
\begin{lemma}\label{lem:rec} Let 
$
(M,d)$ be a metric space. 
If $x',x''\in M$, $x'\neq x''$, $y\in M
\setminus B(x',2d(x',x''))$, then	 
 \[
 \frac{1}{2 }d(x',y)\leq d(x'',y)\leq 2d(x',y)\,,
 \]
\end{lemma}

\section{Composite kernels on upper Ahlfors regular sets}\label{sec:cmpoke}
Next we analyze the composition of kernels on upper Ahlfors regular sets that 
 are continuous outside the diagonal and that have a weak singularity. 
 Here we are extending to the setting of upper Ahlfors regular sets previous results of Giraud \cite{Gi34} (see Miranda \cite[Chap.~II, \S \ 11, p.~24]{Mi70}, Cialdea \cite[\S 1, nr.~4]{Ci01}).
 
 To do so, we introduce the following technical statement, that we prove by exploiting Lemmas \ref{lem:cominur}, \ref{lem:cominu} and an argument that is different from the classical one in case $M={\mathbb{R}}^n$. 
\begin{proposition}\label{prop:xytam} 
Let $X$, $Y$, $Z$ be subsets  of a metric space  $(M,d)$.   Let $\nu$ be as in (\ref{eq:nu}). Let $\nu(Y)<+\infty$. Let $\upsilon_Y\in]0,+\infty[$. 
Let $Y$ be upper $\upsilon_Y$-Ahlfors regular with respect to $X\cup Z$. Let $s_1$, $s_2\in[0,\upsilon_Y[$. If $s_1+s_2=\upsilon_Y$, then we further assume that $Y$ is strongly upper $\upsilon_Y$-Ahlfors regular with respect to $X\cup Z$. 

Then there exists $c_{X,Y,Z}^{(s_1,s_2)}\in]0,+\infty[$ such that
\begin{eqnarray}\label{prop:xytam1}
\lefteqn{
\int_Y\frac{d\nu(y)}{d(x,y)^{s_1}d(y,z)^{s_2}}
}
\\ \nonumber
&& 
\leq\left\{
\begin{array}{ll}
 c_{X,Y,Z}^{(s_1,s_2)}(1+d(x,z)^{\upsilon_Y-(s_1+s_2)}) & 
 \\
 \qquad\qquad\text{if}\ (x,z)\in X\times Z \,,\ s_1+s_2<\upsilon_Y\,,&
 \\
 c_{X,Y,Z}^{(s_1,s_2)}(1+|\log d(x,z)|) & 
 \\
  \qquad\qquad\text{if}\ 
 (x,z)\in (X\times Z)\setminus {\mathbb{D}}_{X\times Z}\,,\ 
 s_1+s_2=\upsilon_Y\,,&
 \\
 c_{X,Y,Z}^{(s_1,s_2)}d(x,z)^{\upsilon_Y-(s_1+s_2)}&
 \\
  \qquad\qquad
  \text{if}\ 
 (x,z)\in (X\times Z)\setminus {\mathbb{D}}_{X\times Z}\,,\ s_1+s_2>\upsilon_Y\,,&
\end{array}
\right.
\end{eqnarray}
for all $(x,z)\in (X\times Z)\setminus {\mathbb{D}}_{X\times Z}$.
\end{proposition}
{\bf Proof.} We first note that if $(x,z)\in (X\times Z)\setminus {\mathbb{D}}_{X\times Z}$, then
\begin{eqnarray*}
\lefteqn{\int_Y\frac{d\nu(y)}{d(x,y)^{s_1}d(y,z)^{s_2}}
}
\\ \nonumber
&&\quad
\leq\int_{B(x,d(x,z)/2)}\frac{d\nu(y)}{d(x,y)^{s_1}d(y,z)^{s_2}}+\int_{B(z,d(x,z)/2)}\frac{d\nu(y)}{d(x,y)^{s_1}d(y,z)^{s_2}}
\\ \nonumber
&&\quad
+\int_{Y\setminus\left(B(x,d(x,z)/2)\cup B(z,d(x,z)/2)
\right)
}\frac{d\nu(y)}{d(x,y)^{s_1}d(y,z)^{s_2}}\,.
\end{eqnarray*}
By the triangular inequality, we have
\begin{eqnarray*}
\lefteqn{
d(y,z)\geq d(z,x)-d(x,y)\geq d(z,x)-d(x,z)/2=d(x,z)/2
}
\\
&&\qquad\qquad\qquad\qquad\qquad\qquad\qquad\qquad
\qquad
\qquad\forall y\in  B(x,d(x,z)/2)\,,
\\ \nonumber
\lefteqn{
d(y,x)\geq d(x,z)-d(z,y)\geq d(x,z)-d(z,x)/2=d(z,x)/2
}
\\
&&\qquad\qquad\qquad\qquad\qquad\qquad\qquad\qquad
\qquad
\qquad\forall y\in  B(z,d(x,z)/2)\,.
\end{eqnarray*}
Then Lemma \ref{lem:cominur} (iii) implies that
\begin{eqnarray}\label{prop:xytam1a}
\lefteqn{\int_Y\frac{d\nu(y)}{d(x,y)^{s_1}d(y,z)^{s_2}}
\leq
\left(
d(x,z)/2
\right)^{-s_2} c_{s_1,X\cup Z,Y}''		\left(
d(x,z)/2
\right)^{\upsilon_Y-s_1}
}
\\ \nonumber
&&\qquad\quad 
+
\left(
d(x,z)/2
\right)^{-s_1} c_{s_2,X\cup Z,Y}''	 	\left(
d(x,z)/2
\right)^{\upsilon_Y-s_2}
\\ \nonumber
&&\qquad\quad
+\int_{\{y\in Y\setminus\left(B(x,d(x,z)/2)\cup B(z,d(x,z)/2)
\right):\,d(y,x)\leq d(y,z)\}
}\frac{d\nu(y)}{d(x,y)^{s_1}d(y,z)^{s_2}}
\\ \nonumber
&&\qquad\quad 
+\int_{\{y\in Y\setminus\left(B(x,d(x,z)/2)\cup B(z,d(x,z)/2)
\right):\,d(y,z)\leq d(y,x)\}
}\frac{d\nu(y)}{d(x,y)^{s_1}d(y,z)^{s_2}}
\\ \nonumber
&&\qquad 
\leq
2\, 2^{s_1+s_2-\upsilon_Y}d(x,z)^{\upsilon_Y-s_1-s_2}
\max\{ c_{s_1,X\cup Z,Y}'' 	
,c_{s_2,X\cup Z,Y}''	 	\}
\\ \nonumber
&&\qquad\quad
+\int_{Y\setminus B(x,d(x,z)/2) }\frac{d\nu(y)}{d(x,y)^{s_1+s_2} }
+\int_{Y\setminus B(z,	d (x,z)/2) }\frac{d\nu(y)}{d(y,z)^{s_1+s_2} }
\,.
\end{eqnarray}
We now consider the three cases of the statement separately. If $s_1+s_2<\upsilon_Y$, then inequality (\ref{prop:xytam1a}) and Lemma \ref{lem:cominur} (ii) imply that
\begin{eqnarray}\label{prop:xytam2}
\lefteqn{\int_Y\frac{d\nu(y)}{d(x,y)^{s_1}d(y,z)^{s_2}}
}
\\ \nonumber
&&\leq
2\,2^{s_1+s_2-\upsilon_Y}d(x,z)^{\upsilon_Y-s_1-s_2}
\max\{ c_{s_1,X\cup Z,Y}''	 
,
c_{s_2,X\cup Z,Y}''	 	\} 
\\ \nonumber
&&\quad
+2c_{s_1+s_2,X\cup Z,Y}' \qquad\forall (x,z)\in (X\times Z)\setminus {\mathbb{D}}_{X\times Z}\,.
\end{eqnarray}
On the other hand, if 
$(x,z)\in {\mathbb{D}}_{X\times Z}$, then Lemma \ref{lem:cominur} (ii) implies that
\begin{equation}\label{prop:xytam2a}
\int_Y\frac{d\nu(y)}{d(x,y)^{s_1}d(y,z)^{s_2}}
=\int_Y\frac{d\nu(y)}{d(x,y)^{s_1+s_2} }
\leq c'_{s_1+s_2,X,Y}\,.
\end{equation}
If $s_1+s_2=\upsilon_Y$ and if $(x,z)\in (X\times Z)\setminus {\mathbb{D}}_{X\times Z}$,  $d(x,z)\geq 1/e$, then inequality (\ref{prop:xytam1a}) implies that
\begin{equation}\label{prop:xytam3}
 \int_Y\frac{d\nu(y)}{d(x,y)^{s_1}d(y,z)^{s_2}}
\leq
2\max\{ c_{s_1,X\cup Z,Y}''	 	
,
 c_{s_2,X\cup Z,Y}''	 	\} 
+2\nu(Y)(1/(2e))^{-\upsilon_Y}\,.
\end{equation}
If instead $(x,z)\in (X\times Z)\setminus {\mathbb{D}}_{X\times Z}$, $d(x,z)< 1/e$, then inequality (\ref{prop:xytam1a}) and Lemma \ref{lem:cominu} (ii) imply that
\begin{eqnarray}\label{prop:xytam4}
\lefteqn{
 \int_Y\frac{d\nu(y)}{d(x,y)^{s_1}d(y,z)^{s_2}}
 }
 \\ \nonumber
 &&\qquad
\leq
2\max\{ c_{s_1,X\cup Z,Y}''	 	
,
 c_{s_2,X\cup Z,Y}''	 	\} 
+
2c^{(iv)}_{X\cup Z,Y}|\log(d(x,z)/2)|
\,.
\end{eqnarray}
If $s_1+s_2>\upsilon_Y$, then inequality (\ref{prop:xytam1a}) and Lemma \ref{lem:cominu} (i) imply that
\begin{eqnarray}\label{prop:xytam5}
\lefteqn{\int_Y\frac{d\nu(y)}{d(x,y)^{s_1}d(y,z)^{s_2}}
}
\\ \nonumber
&& 
\leq 2\,
2^{s_1+s_2-\upsilon_Y}d(x,z)^{\upsilon_Y-s_1-s_2}
\max\{ c_{s_1,X\cup Z,Y}'' 
,
 c_{s_2,X\cup Z,Y}''	 \} 
\\ \nonumber
&&\quad 
+2 c'''_{s_1+s_2,X\cup Z,Y}d(x,z)^{\upsilon_Y-s_1-s_2} 2^{s_1+s_2-\upsilon_Y}
  \end{eqnarray}
for all $ (x,z)\in (X\times Z)\setminus {\mathbb{D}}_{X\times Z}$. 
Then by combining inequalities (\ref{prop:xytam2})--(\ref{prop:xytam5}), we conclude that the constant $c_{X,Y,Z}^{(s_1,s_2)}$ of inequality (\ref{prop:xytam1}) does exist.\hfill  $\Box$ 

\vspace{\baselineskip}

Thus we are ready to prove the following extension of a classical result on composite kernels (see also Miranda \cite[Chap II, \S 11, p.~24]{Mi70}).
\begin{theorem}
\label{thm:cokewa}
Let  $Y$  be a subset  of a metric space  $(M,d)$.   Let $\nu$ be as in (\ref{eq:nu}) with $X=Y$. Let $\nu(Y)<+\infty$.   Let $\upsilon_Y\in]0,+\infty[$. 
Let $Y$ be upper $\upsilon_Y$-Ahlfors regular with respect to $Y$. Let $s_1$, $s_2\in[0,\upsilon_Y[$. If $s_1+s_2=\upsilon_Y$, then we further assume that $Y$ is strongly upper $\upsilon_Y$-Ahlfors regular with respect to $Y$. 

Let $K_{l} \in {\mathcal{K}}_{s_l, Y\times Y}$ for $l\in\{1,2\}$.  Then the following statements hold
 \begin{enumerate}
\item[(i)] If $s_{1}+s_{2}\geq \upsilon_Y$, then for each $(x,y)\in Y\times Y\setminus {\mathbb{D}}_{   Y\times Y  }$, the function from  $Y$ to ${\mathbb{C}}$ which takes $t$ to $ K_{1}(x,t)K_{2}(t,y)$ is $\nu$-integrable in $Y$ and the function $K_{3}$ from $Y^{2}\setminus {\mathbb{D}}_{ Y\times Y 	 }$ to ${\mathbb{C}}$ defined by
\begin{equation}\label{thm:cokewa0}
K_{3}(x,y)\equiv\int_{Y}K_{1}(x,t)K_{2}(t,y)\,d\nu(t)
\qquad\forall (x,y)\in Y^{2}\setminus {\mathbb{D}}_{ Y\times Y  }\,,
\end{equation}
 is continuous on $Y^{2}\setminus {\mathbb{D}}_{    Y\times Y   }$. Moreover,   
\begin{eqnarray}\label{thm:cokewa1}
\lefteqn{
\sup_{(x,y)\in 
Y^{2}\setminus {\mathbb{D}}_{   Y\times Y   }}|K_{3}(x,y)|d(x,y)^{s_{1}+s_{2}-\upsilon_Y}
 <+\infty
 }
 \\ \nonumber
&&\qquad\qquad\qquad\qquad\qquad\qquad\qquad {\text{if}}\ s_{1}+s_{2}> \upsilon_Y\,,
\\ \nonumber
\lefteqn{
\sup_{(x,y)\in 
Y^{2}\setminus {\mathbb{D}}_{   Y\times Y   }} |K_{3}(x,y)|(1+|\ln d(x,y)|)^{-1}<+\infty
}
\\ \nonumber
&&\qquad\qquad\qquad\qquad\qquad\qquad\qquad {\text{if}}\ s_{1}+s_{2}= \upsilon_Y\,.
\end{eqnarray}
  
\item[(ii)]  If $s_{1}+s_{2}< \upsilon_Y$, then for each $(x,y)\in 
Y^{2}$, the function from $Y$ to ${\mathbb{C}}$ which takes $t$ to $ K_{1}(x,t)K_{2}(t,y)$ is $\nu$-integrable in $Y$, and the function $\tilde{K}_{3}$ from 
$Y^{2}$ to ${\mathbb{C}}$ defined by
\[
\tilde{K}_{3}(x,y)\equiv\int_{Y}K_{1}(x,t)K_{2}(t,y)\,d\nu(t)
\qquad\forall (x,y)\in Y^{2}\,,
\]
is continuous. Moreover,   
\begin{equation}\label{thm:cokewa2}
\sup_{(x,y)\in 
Y^{2} }|\tilde{K}_{3}(x,y)|
(1+d(x,y)^{\upsilon_Y-(s_{1}+s_{2})})^{-1}<+\infty\,.
\end{equation}
\end{enumerate}
\end{theorem}
{\bf Proof.} Since both $K_{1}$ and $K_{2}$  are continuous, then  $K_{1}(x,\cdot)
K_{2}(\cdot,y)$ is continuous in $Y\setminus\{x,y\}$ for all $x,y\in Y$. By assumption, we have
\[
|K_{1}(x,t)K_{2}(t,y)|\leq
\frac{\|K_1\|_{ {\mathcal{K}}_{s_1, Y\times Y} }
}{  d(x,t)^{s_{1}} }\frac{\|K_2\|_{ {\mathcal{K}}_{s_2, Y\times Y} }
}{  d(y,t)^{s_{2}} }
\qquad\forall t\in Y\setminus\{x,y\}\,,
\]
for all $(x,y)\in Y^{2}$. Then Proposition \ref{prop:xytam} implies the integrability of the function $K_{1}(x,\cdot)
K_{2}(\cdot,y)$ in $Y$ for all $(x,y)\in Y^{2}\setminus {\mathbb{D}}_{ Y\times Y   }$ if $s_{1}+s_{2}\geq\upsilon_Y$ and for all
$(x,y)\in Y^{2}$ if $s_{1}+s_{2}<\upsilon_Y$. Proposition \ref{prop:xytam} implies also the validity of inequalities (\ref{thm:cokewa1}), (\ref{thm:cokewa2}). Next we turn to the problem of continuity of $K_{3}$, $\tilde{K}_{3}$. 

We first consider case $s_{1}+s_{2}\geq \upsilon_Y$.
Let $(\tilde{x},\tilde{y})\in Y^{2}$, $\tilde{x}\neq \tilde{y}$. We want to prove that
\[
\lim_{(x,y)\to (\tilde{x},\tilde{y}) }|K_{3}(x,y)-K_{3}(\tilde{x},\tilde{y})|=0\,.
\]
By Lemma  \ref{lem:cominur} (ii)  there exists $\rho_{\epsilon}\in ]0,d(\tilde{x},\tilde{y})/4[$ such that
\[
\sup_{x\in Y}\int_{Y
\cap B(\tilde{x},\rho_{\epsilon})}
\frac{1}{d(x,t)^{s_{1}}}\,d\nu(t)\leq\epsilon\,,
\quad
\sup_{y\in Y }\int_{Y
\cap B(\tilde{y},\rho_{\epsilon})}
\frac{1}{d(y,t)^{s_{2}}}\,d\nu(t)\leq\epsilon 
\,.
\]
Clearly,
\begin{eqnarray*}
\lefteqn{
|K_{3}(x,y)-K_{3}(\tilde{x},\tilde{y})|
}
\\
\nonumber
&&
=\left|
\int_{Y}
K_{1}(x,t)K_{2}(t,y)\,d\nu(t)
-
\int_{Y}
K_{1}(\tilde{x},t)K_{2}(t,\tilde{y})\,d\nu(t)
\right|
\leq
\sum_{l=1}^{3}I_{l}(x,y)\,,
\end{eqnarray*}
for all $x,y\in Y$, $x\neq y$, where
\begin{eqnarray*}
I_{1}(x,y)&\equiv &
\int_{Y
\setminus( B(\tilde{x},\rho_{\epsilon})
\cup B(\tilde{y},\rho_{\epsilon}))
}
|K_{1}(x,t)K_{2}(t,y)-K_{1}(\tilde{x},t)K_{2}(t,\tilde{y})|\,d\nu(t)\,,
\\
I_{2}(x,y)&\equiv &
\int_{Y
\cap B(\tilde{x},\rho_{\epsilon})
}
|K_{1}(x,t)K_{2}(t,y)-K_{1}(\tilde{x},t)K_{2}(t,\tilde{y})|\,d\nu(t)\,,
\\
I_{3}(x,y)&\equiv &
\int_{Y
\cap  B(\tilde{y},\rho_{\epsilon})
}
|K_{1}(x,t)K_{2}(t,y)-K_{1}(\tilde{x},t)K_{2}(t,\tilde{y})|\,d\nu(t)\,,
\end{eqnarray*}
for all $(x,y)\in Y^{2}\setminus {\mathbb{D}}_{ Y\times Y}$. Since $K_1$ and $K_2$ are continuous outside of the diagonal,
\[
|K_{1}(x,t)K_{2}(t,y)-K_{1}(\tilde{x},t)K_{2}(t,\tilde{y})|
\leq 2  \frac{
\|K_1\|_{ {\mathcal{K}}_{s_1, Y\times Y} }\|K_2\|_{ {\mathcal{K}}_{s_2, Y\times Y} }}
{(\rho_{\epsilon}/2)^{s_1+s_2} }
\]
for all 
\[
(x,y,t)\in (Y\cap B(\tilde{x},\rho_{\epsilon}/2))
\times
(Y\cap B(\tilde{y},\rho_{\epsilon}/2))
\times (  Y
\setminus( B(\tilde{x},\rho_{\epsilon})
\cup B(\tilde{y},\rho_{\epsilon})) )\,,
\] 
the function in the argument of the integral in 
$I_{1}(x,y)$ is continuous and bounded in 
\[
(x,y,t)\in (Y\cap B(\tilde{x},\rho_{\epsilon}/2))
\times
(Y\cap B(\tilde{y},\rho_{\epsilon}/2))
\times (  Y
\setminus( B(\tilde{x},\rho_{\epsilon})
\cup B(\tilde{y},\rho_{\epsilon})) )\,.
\]
and $Y$ has finite measure, then the classical theorem of continuity for integrals depending on a parameter implies that
\[
\lim_{(x,y)\to (\tilde{x},\tilde{y}) }I_{1}(x,y)=0\,.
\]
We now consider $I_{2}(x,y)$. By assumption, we have
\begin{eqnarray*}
\lefteqn{
|I_{2}(x,y)|\leq \|K_1\|_{ {\mathcal{K}}_{s_1, Y\times Y} }\|K_2\|_{ {\mathcal{K}}_{s_2, Y\times Y} }
}
\\ \nonumber
&&\qquad\qquad\times
\int_{Y
\cap B(\tilde{x},\rho_{\epsilon})
}
\left(\frac{1}{d(x,t)^{s_{1}}  d(y,t)^{s_{2}}   }
+
\frac{1}{d(\tilde{x},t)^{s_{1}}   d(\tilde{y},t)^{s_{2}}   }\right)\,d\nu(t)\,,
\end{eqnarray*}
for all $x,y\in Y$. Next we note that if $d(y,\tilde{y})\leq \rho_{\epsilon}$, then
\[
d(\tilde{x},\tilde{y})\leq d(\tilde{x},t)+d(t,y)+d(y,\tilde{y})
\leq 2 \rho_{\epsilon}+d(t,y)
\leq
\frac{1}{2}d(\tilde{x},\tilde{y})+d(t,y)\,,
\]
for all $t\in  B(\tilde{x},\rho_{\epsilon})$, and accordingly
\[
d(y,t)>\frac{1}{2}d(\tilde{x},\tilde{y})\qquad\forall t\in B(\tilde{x},\rho_{\epsilon})\,.
\] 
Hence, our choice of $\rho_{\epsilon}$ implies that
\begin{eqnarray*}
\lefteqn{|I_{2}(x,y)| 
\leq \|K_1\|_{ {\mathcal{K}}_{s_1, Y\times Y} }\|K_2\|_{ {\mathcal{K}}_{s_2, Y\times Y} }\left(\frac{2}{d(\tilde{x},\tilde{y})}\right)^{s_{2}}
}
\\
\nonumber
&&\qquad\qquad\qquad\quad\times
\int_{Y
\cap  B(\tilde{x},\rho_{\epsilon})
}
\left(\frac{1}{d(x,t)^{s_{1}}      }
+
\frac{1}{d(\tilde{x},t)^{s_{1}}     }\right)\,d\nu(t)
\\
\nonumber
&&\qquad\ \ 
\leq
2 \|K_1\|_{ {\mathcal{K}}_{s_1, Y\times Y} }\|K_2\|_{ {\mathcal{K}}_{s_2, Y\times Y} }
\left(\frac{2}{d(\tilde{x},\tilde{y})}\right)^{s_{2}}\epsilon\,,
\end{eqnarray*}
for all $(x,y)\in Y^{2}$ such that $x\neq y$, $d(y,\tilde{y})\leq\rho_{\epsilon}$.
By interchanging the roles of $\tilde{x}$ and $\tilde{y}$, we can prove that
\[
|I_{3}(x,y)| \leq
2 \|K_1\|_{ {\mathcal{K}}_{s_1, Y\times Y} }\|K_2\|_{ {\mathcal{K}}_{s_2, Y\times Y} }\left(\frac{2}{d(\tilde{x},\tilde{y})}\right)^{s_{1}}\epsilon\,,
\]
for all $(x,y)\in Y^{2}$ such that $x\neq y$, $d(x,\tilde{x})\leq\rho_{\epsilon}$. Hence,
\begin{eqnarray*}
\lefteqn{
\limsup_{(x,y)\to (\tilde{x},\tilde{y})}|K_{3}(x,y)-K_{3}(\tilde{x},\tilde{y})|
\leq 
\limsup_{(x,y)\to (\tilde{x},\tilde{y})}|I_{1}(x,y) |}
\\
\nonumber
&&
+
2 \|K_1\|_{ {\mathcal{K}}_{s_1, Y\times Y} }\|K_2\|_{ {\mathcal{K}}_{s_2, Y\times Y} }
\left[
\left(\frac{2}{d(\tilde{x},\tilde{y})}\right)^{s_{2}}
+
\left(\frac{2}{d(\tilde{x},\tilde{y})}\right)^{s_{1}}
\right]\epsilon\,.
\end{eqnarray*}
Since $\epsilon>0$ is arbitrary, we have 
$\limsup_{(x,y)\to (\tilde{x},\tilde{y})}|K_{3}(x,y)-K_{3}(\tilde{x},\tilde{y})|=0$ and accordingly
$\lim_{(x,y)\to (\tilde{x},\tilde{y})}|K_{3}(x,y)-K_{3}(\tilde{x},\tilde{y})|$=0, and the proof of case $s_{1}+s_{2}\geq \upsilon_Y$ is complete.\par

We now consider case $s_{1}+s_{2}< \upsilon_Y$. Let $(\tilde{x},\tilde{y})\in Y^{2}$.  We want to prove that
\[
\lim_{(x,y)\to (\tilde{x},\tilde{y}) }|\tilde{K}_{3}(x,y)-\tilde{K}_{3}(\tilde{x},\tilde{y})|=0\,.
\]
To do so, we plan to apply the Vitali Convergence Theorem in $Y$ and to prove that if $\{(x_{j},y_{j})\}_{j\in {\mathbb{N}} }$ is a sequence in  $Y^{2}\setminus\{(\tilde{x},\tilde{y})\}$, then
\begin{equation}
\label{thm:cokewa2a}
\lim_{j\to \infty}
\int_{Y} K_{1}(x_{j},t)K_{2}(t,y_{j})\,d\nu(t)
=
\int_{Y} K_{1}(\tilde{x},t)K_{2}(t,\tilde{y})\,d\nu(t)\,,
\end{equation}
(cf.,~e.g., Folland~\cite[Ex.~15, p.~187]{Fo99}, Fichera~\cite{Fi43}). By assumption, if $t\in Y\setminus\{\tilde{x},\tilde{y}\}$, then 
$K_{1}(x,t)K_{2}(t,y)$ is continuous when $(x,y)$ ranges in a neighborhood of $(\tilde{x},\tilde{y})$ and accordingly
\[
\lim_{j\to \infty}
 K_{1}(x_{j},t)K_{2}(t,y_{j})=
  K_{1}(\tilde{x},t)K_{2}(t,\tilde{y}) \,.
\]
In particular,  $\{ K_{1}(x_{j},t)K_{2}(t,y_{j})\}_{j\in {\mathbb{N}} }$ converges  to $K_{1}(\tilde{x},t)K_{2}(t,\tilde{y})$
for $\nu$-almost all $t\in Y$ (cf.~Lemma \ref{lem:cominur} (i)). Then it suffices to show that for each $\epsilon>0$, there exists $\delta>0$ such that
\[
\int_{E}\left| K_{1}(x_{j},t)K_{2}(t,y_{j})
\right|
\,d\nu(t) 
\leq \epsilon\qquad\forall j\in {\mathbb{N}}\,,
\]
for all   $E\in {\mathcal{N}}$ such that $\nu(E)\leq \delta$. To do so, we look for 
$p,q\in ]1,+\infty[$ such that
\begin{equation}
\label{thm:cokewa3}
\frac{1}{p}+\frac{1}{q}=1\,,\qquad
ps_{1}<\upsilon_Y\,,
\qquad
qs_{2}<\upsilon_Y\,.
\end{equation}
Since $\lim_{ p\to (\upsilon_Y/s_{1})^{-} }1-(1/p)=\frac{\upsilon_Y-s_{1}}{\upsilon_Y}
>\frac{s_{2}}{\upsilon_Y}$, then there exists $p\in ]1,\upsilon_Y/s_{1}[$ such that $1-(1/p)> s_{2}/\upsilon_Y$. Then we have 
$
ps_{1}<\upsilon_Y$ and  
\[
qs_{2}=\frac{p}{p-1} s_{2}=\frac{1}{1-(1/p)}s_{2}<\frac{1}{(s_{2}/\upsilon_Y)}s_{2}=\upsilon_Y
\]
 and (\ref{thm:cokewa3}) holds true. Then we choose
$p$, $q$ as in (\ref{thm:cokewa3}). By the H\"{o}lder inequality, we have
\begin{eqnarray*}
\lefteqn{
\int_{E} \left|K_{1}(x_{j},t)K_{2}(t,y_{j})\right|\,d\nu(t)
}
\\
\nonumber
&&
\leq
\|K_1\|_{ {\mathcal{K}}_{s_1, Y\times Y} }\|K_2\|_{ {\mathcal{K}}_{s_2, Y\times Y} }
\int_{E}
\frac{1}{   d(x_{j},t)^{s_{1}}   }
\frac{1}{d(y_{j},t)^{s_{2}}       }
\,d\nu(t)
\\
\nonumber
&& 
\leq \|K_1\|_{ {\mathcal{K}}_{s_1, Y\times Y} }\|K_2\|_{ {\mathcal{K}}_{s_2, Y\times Y} }
\\
\nonumber
&&\quad\times 
\left(
\int_{E}
\frac{1}{  d(x_{j},t)^{ps_{1}}   }\,d\nu(t)
\right)^{1/p}
\left(
\int_{E}
\frac{1}{   d(y_{j},t)^{qs_{2}}   }\,d\nu(t)
\right)^{1/q}
\end{eqnarray*}
for all $j\in {\mathbb{N}}$.
Since $s_{1} p<\upsilon_Y$,  $s_{2} q<\upsilon_Y$, then there exists $\delta\in]0,+\infty[$ such that
\[
\int_{E}
\frac{1}{  d(x_{j},t)^{ps_{1}}   }\,d\nu(t)\leq\epsilon\,,
\qquad
\int_{E}
\frac{1}{   d(y_{j},t)^{qs_{2}}   }\,d\nu(t)\leq \epsilon
\qquad
\forall j\in {\mathbb{N}}\,,
\]
for all  $E\in {\mathcal{N}}$  such that $\nu(E)\leq \delta$ (see Lemma \ref{lem:cominuu}). Then we also have
\begin{eqnarray*}
\lefteqn{
\int_{E} \left|K_{1}(x_{j},t)K_{2}(t,y_{j}) \right| \,d\nu(t)
}
\\ \nonumber
&&\qquad
\leq \|K_1\|_{ {\mathcal{K}}_{s_1, Y\times Y} }\|K_2\|_{ {\mathcal{K}}_{s_2, Y\times Y} } \epsilon^{1/p}\epsilon^{1/q}= \|K_1\|_{ {\mathcal{K}}_{s_1, Y\times Y} }\|K_2\|_{ {\mathcal{K}}_{s_2, Y\times Y} }\epsilon\,,
\end{eqnarray*}
for all  $E\in{\mathcal{N}}$ such that $\nu(E)\leq \delta$. Then the Vitali Convergence Theorem implies that 
 (\ref{thm:cokewa2a})  holds and that $\tilde{K}_{3}$ is continuous on $Y^{2}$.  \hfill  $\Box$ 

\vspace{\baselineskip}

 In order to get a better understanding of the order of singularity of composite kernels, we introduce the following immediate consequence  of statement (i) of the previous Theorem \ref{thm:cokewa}.
\begin{corollary}\label{corol:cokewa}
 Let  $Y$  be a subset  of a metric space  $(M,d)$.    Let $\nu$ be as in (\ref{eq:nu}) with $X=Y$. Let $\nu(Y)<+\infty$.   Let $\upsilon_Y\in]0,+\infty[$.
Let $Y$ be upper $\upsilon_Y$-Ahlfors regular with respect to $Y$. Let $t_1$, $t_2\in]0,\upsilon_Y]$. Let $K_{l} \in {\mathcal{K}}_{\upsilon_Y-t_{l}, Y\times Y}$ for $l\in\{1,2\}$. 
 If $t_{1}+t_{2}<\upsilon_Y$, then  the composite kernel $K_{3}$ of $K_{1}$ and $K_{2}$ is continuous in $Y^{2}\setminus {\mathbb{D}}_{   Y\times Y   }$ and satisfies the following inequality
\[
\sup_{(x,y)\in 
Y^{2}\setminus {\mathbb{D}}_{    Y\times Y   }}|K_{3}(x,y)|d(x,y)^{\upsilon_Y-(t_{1}+t_{2})}
 <+\infty\,.
\]
\end{corollary}
{\bf Proof.} It suffices to set $s_{l}\equiv \upsilon_Y-t_{l}$
for $l=1,2$ and to invoke statement (i) of  Theorem \ref{thm:cokewa}.\hfill  $\Box$ 

\vspace{\baselineskip}

Then by the previous Theorem \ref{thm:cokewa} and Corollary \ref{corol:cokewa}, we readily deduce the validity of the following statement.
\begin{corollary}\label{corol:iterwsa}
 Let  $Y$  be a bounded subset  of a metric space  $(M,d)$.   Let $\nu$ be as in (\ref{eq:nu}) with $X=Y$.    Let $\upsilon_Y\in]0,+\infty[$. Let $s\in [0,\upsilon_Y[$.  
Let $Y$ be upper $\upsilon_Y$-Ahlfors regular with respect to $Y$. 
Let $K\in {\mathcal{K}}_{s,Y\times Y}$.
  
Then there exists $r\in {\mathbb{N}}\setminus\{0\}$ such that $K^{(r)}$ has a continuous and bounded extension to  $Y^{2}$. 
\end{corollary}
{\bf Proof.}   By Remark \ref{rem:bddimfi} with $X=Y$, we have $\nu(Y)<+\infty$. 
 If $s=0$, then  Theorem \ref{thm:cokewa} (ii) implies that $K^{(2)}$ admits a continuous and bounded extension to $Y^2$. Thus we can assume that $s\in ]0,\upsilon_Y[$.	 	
By assumption,   
\[
|K(x,y)|\leq\frac{
\|K\|_{{\mathcal{K}}_{s,Y\times Y}}
}{d(x,y)^{s}}
=\frac{\|K\|_{{\mathcal{K}}_{s,Y\times Y}}}{d(x,y)^{\upsilon_Y-[\upsilon_Y-s]}}
\qquad\forall (x,y)\in Y^{2}\setminus {\mathbb{D}}_{  Y\times Y  }\,.
\]
If $2[\upsilon_Y-s]<\upsilon_Y$, \textit{i.e.}, if $\upsilon_Y<2s$, then the previous Corollary \ref{corol:cokewa}  implies the continuity of $K^{(2)}$ in $ Y^{2}\setminus {\mathbb{D}}_{    Y\times Y}   	 $ and the existence of $d_{2}>0$ such that
\[
|K^{(2)}(x,y)|\leq \frac{d_{2}}{ d(x,y)^{\upsilon_Y-2[\upsilon_Y-s]
}}\qquad\forall (x,y)\in 
Y^{2}\setminus {\mathbb{D}}_{ Y\times Y }\,.
\]
 More generally, if	  $m\in {\mathbb{N}}\setminus\{0\}$ and   $m
[\upsilon_Y-s]<\upsilon_Y$, then	  a simple inductive argument based on the previous Corollary \ref{corol:cokewa} shows that   $K^{(m)}$ is continuous in $ Y^{2}\setminus {\mathbb{D}}_{    Y\times Y}   	 $ and that there exists $d_{m}>0$ such that
\[
|K^{(m)}(x,y)|\leq \frac{d_{m}}{ d(x,y)^{\upsilon_Y-m[\upsilon_Y-s]
}}\qquad\forall (x,y)\in 
Y^{2}\setminus {\mathbb{D}}_{ Y\times Y }\,.
\]
Now let $\tilde{m}$ be the maximum $m\in {\mathbb{N}}\setminus\{0\}$ such that
\[
m
[\upsilon_Y-s]<\upsilon_Y\,.
\]
Clearly,
\[
(\tilde{m}+1)[\upsilon_Y-s]\geq \upsilon_Y\,.
\]
If the strict inequality holds, then we have $\tilde{m}[\upsilon_Y-s]>s$ and accordingly
\[
s +\{
\upsilon_Y-\tilde{m}[\upsilon_Y-s]
\}<s+\upsilon_Y-s=\upsilon_Y\,.
\]
Hence,  the assumption that $Y$ is bounded and  	Theorem \ref{thm:cokewa} (ii)  imply  	that $K^{(\tilde{m}+1)}$ admits a continuous and bounded extension to $Y^2$.  

If instead $(\tilde{m}+1)[\upsilon_Y-s]= \upsilon_Y$, we choose $s'\in ]s,\upsilon_Y[$. Since $Y$ is bounded, we have
\[
K\in {\mathcal{K}}_{\upsilon_Y-[\upsilon_Y-s'],Y\times Y}\,.
\]
Then the membership of $K^{(\tilde{m})}$ in ${\mathcal{K}}_{\upsilon_Y-\tilde{m}[\upsilon_Y-s],Y\times Y}$, the inequality
\[
\tilde{m}[\upsilon_Y-s]+[\upsilon_Y-s']=\tilde{m}[\upsilon_Y-s]+[\upsilon_Y-s]+(s-s')=\upsilon_Y+(s-s')<\upsilon_Y
\]
and Corollary \ref{corol:cokewa} imply that
\[
K^{(\tilde{m}+1)}\in {\mathcal{K}}_{\upsilon_Y-\tilde{m}[\upsilon_Y-s]-[\upsilon_Y-s'],Y\times Y}\,.
\]
Next, we note that
\[
s+\{
\upsilon_Y-\tilde{m}[\upsilon_Y-s]-[\upsilon_Y-s']
\}
=\upsilon_Y-\tilde{m}[\upsilon_Y-s]-[\upsilon_Y-s] +s'=0+s'<\upsilon_Y	 \,.	 
\]
Hence,  the assumption that $Y$ is bounded and  	Theorem \ref{thm:cokewa} (ii)  imply  	that  $K^{(\tilde{m}+1)+1}$ admits a continuous and bounded extension to $Y^2$.\hfill  $\Box$ 

\vspace{\baselineskip}

Next we prove the following technical statement.
\begin{lemma}\label{lem:sxytam}
 Let  $Y$  be a subset  of a metric space  $(M,d)$.   Let $\nu$ be as in (\ref{eq:nu}) with $X=Y$.   Let $\upsilon_Y\in]0,+\infty[$.  Let  $Y$ be  upper $\upsilon_Y$-Ahlfors regular with respect to $Y$. Let  $\nu(Y)<+\infty$ whenever $r_{Y,Y,\upsilon_Y}<+\infty$. Let $s_1$, $s_2\in[0,\upsilon_Y[$, $a\in]0,+\infty[$. Then there exists $c_{s_1,s_2,a,Y}\in]0,+\infty[$ such that
 \begin{eqnarray}\label{lem:sxytam1}
\lefteqn{
\int_{Y\cap B(\xi,ad(x',x''))}\frac{d\nu(\eta)}{d(\xi,\eta)^{s_1}d(\eta,y)^{s_2}}
}
\\ \nonumber
&&\qquad
\leq c_{s_1,s_2,a,Y}\left(d(	x'	,y)^{-s_2}d(x',x'')^{\upsilon_Y-s_1}+d(x',y)^{-s_1}d(x',x'')^{\upsilon_Y-s_2}\right)
\end{eqnarray}
 for all $x'$, $x''\in Y$, $x'\neq x''$, $\xi\in \{x',x''\}$, $y\in Y\setminus B(x',2d(x',x''))$. 
 \end{lemma}
{\bf Proof.} We first observe that
\begin{equation}\label{lem:sxytam2}
\max\{d(y_1,\eta),d(\eta,y_2)\}\geq d(y_1,y_2)/2\qquad\forall y_1,y_2,\eta\in Y\,.
\end{equation}
Indeed, if there exist $y_1,y_2,\eta\in Y$ such that inequality (\ref{lem:sxytam2}) is violated, then
\[
d(y_1,y_2)\leq d(y_1,\eta)+d(\eta,y_2)<(d(y_1,y_2)/2)+(d(y_1,y_2)/2)=d(y_1,y_2)\,,
\]
a contradiction. Next we observe that
\begin{equation}\label{lem:sxytam3}
B(\xi,d(x',x'')/2)\cap B(y,d(x',x'')/2)=\emptyset
\end{equation}
for all $x'$, $x''\in Y$, $x'\neq x''$, $\xi\in \{x',x''\}$, $y\in Y\setminus B(x',2d(x',x''))$. Indeed, it there exists $\eta\in B(\xi,d(x',x'')/2)\cap B(y,d(x',x'')/2)$, then
\[
d(\xi,y)\leq d(\xi,\eta)+d(\eta,y)<(d(x',x'')/2)+(d(x',x'')/2)=d(x',x'')\,.
\]
If $\xi=x'$, then our assumption that $d(x',y)\geq 2d(x',x'')$ yields a contradiction. If instead 
$\xi=x''$, then the triangular inequality implies that
\[
2d(x',x'')\leq d(x',y)\leq d(x',x'')+d(x'',y)=d(x',x'')+d(\xi,y)
\]
and accordingly that $d(x',x'')\leq d(\xi,y)<d(x',x'')$ and thus again a contradiction. By inequality (\ref{lem:sxytam2}) and equality (\ref{lem:sxytam3}), we obtain that
\begin{eqnarray}\label{lem:sxytam4}
\lefteqn{
\int_{Y\cap B(\xi,ad(x',x''))}\frac{d\nu(\eta)}{d(\xi,\eta)^{s_1}d(\eta,y)^{s_2}}
}
\\ \nonumber
&& 
\leq
\int_{Y\cap B(\xi,d(x',x'')/8)}\frac{d\nu(\eta)}{d(\xi,\eta)^{s_1}d(\eta,y)^{s_2}}
\\ \nonumber
&& \quad+\int_{Y\cap B(y,d(x',x'')/8)}\frac{d\nu(\eta)}{d(\xi,\eta)^{s_1}d(\eta,y)^{s_2}}
\\ \nonumber
&& \quad+\int_{A(x',x'',a)\cap\{\eta\in Y: d(\xi,\eta)\geq d(\xi,y)/2\}}\frac{d\nu(\eta)}{d(\xi,\eta)^{s_1}d(\eta,y)^{s_2}}
\\ \nonumber
&& \quad+\int_{A(x',x'',a)\cap\{\eta\in Y: d(y,\eta)\geq d(\xi,y)/2\}}\frac{d\nu(\eta)}{d(\xi,\eta)^{s_1}d(\eta,y)^{s_2}}\,,
\end{eqnarray}
where
\[
A(x',x'',a)\equiv  B(\xi,ad(x',x''))\setminus\left(
B(\xi,d(x',x'')/8)\cup B(y,d(x',x'')/8)
\right)
\]
for all $x'$, $x''\in Y$, $x'\neq x''$, $\xi\in \{x',x''\}$, $y\in Y\setminus B(x',2d(x',x''))$. 
Next we turn to estimate the first integral in the right hand side of (\ref{lem:sxytam4}). Since
\[
d(\xi,\eta)<d(x',x'')/8\leq d(x',y)/ 16	 	\leq d(x'',y)/8
\]
(cf.~Lemma \ref{lem:rec}), we have
\begin{eqnarray*}
&&d(x',y)\leq d(x',\eta)+d(\eta,y)\leq(d(x',y)/ 16	)+d(\eta,y)\,,
\\
&&d(x'',y)\leq d(x'',\eta)+d(\eta,y)\leq(d(x'',y)/8	)+d(\eta,y)\,,
\end{eqnarray*}
and accordingly,
\[
 15	 d(x',y)/ 16	 \leq d(\eta,y)\,,\qquad  7	 d(x'',y)/ 8	\leq d(\eta,y)\,,
\]
and thus
\[
d(\xi,y)
 \leq\max\{16/15,8/7\}	d(\eta,y) 
\leq 2 d(\eta,y)
\]
for all $x'$, $x''\in Y$, $x'\neq x''$, $\xi\in \{x',x''\}$, $y\in Y\setminus B(x',2d(x',x''))$, $\eta\in Y\cap B(\xi,d(x',x'')/8)$. Hence, Lemma \ref{lem:cominur} (iii) implies that
\begin{eqnarray}\label{lem:sxytam5}
\lefteqn{
\int_{Y\cap B(\xi,d(x',x'')/8)}\frac{d\nu(\eta)}{d(\xi,\eta)^{s_1}d(\eta,y)^{s_2}}
}
\\ \nonumber
&&\qquad
\leq 2^{s_2} d(\xi,y)^{-s_2}\int_{Y\cap B(\xi,d(x',x'')/8)}\frac{d\nu(\eta)}{d(\xi,\eta)^{s_1}}
\\ \nonumber
&&\qquad
\leq 2^{s_2} d(\xi,y)^{-s_2} c''_{s_1,Y,Y}8^{-(\upsilon_Y-s_1)}d(x',x'')^{\upsilon_Y-s_1}
\,,
\end{eqnarray}
for all $x'$, $x''\in Y$, $x'\neq x''$, $\xi\in \{x',x''\}$, $y\in Y\setminus B(x',2d(x',x''))$.  Next we turn to estimate the second integral in the right hand side of (\ref{lem:sxytam4}). Since
\[
d(y,\eta)<d(x',x'')/8\leq d(x',y)/ 16\leq d(x'',y)/8
\]
(cf.~Lemma \ref{lem:rec}), we have
\begin{eqnarray*}
\lefteqn{d(x',y)\leq d(x',\eta)+d(\eta,y)
}
\\ \nonumber
&&\qquad
\leq
d(x',\eta)+(d(x',x'')/8)\leq 
d(x',\eta)+(d(x',y)/ 16	)
\,,
\\
\lefteqn{d(x'',y)\leq d(x'',\eta)+d(\eta,y)
}
\\ \nonumber
&&\qquad
\leq
d(x'',\eta)+(d(x',x'')/8)\leq 
d(x'',\eta)+(d(x'',y)/ 8	)\,,
\end{eqnarray*}
and accordingly,
\[
\frac{15}{16}	 	d(x',y)\leq d(x',\eta)\,,
\qquad
 \frac{7}{8}	d(x'',y)\leq d(x'',\eta)\,,
\]
and thus
\[
 d(\xi,y) 
  \leq\max\{16/15,8/7\}	d(\xi,\eta)
 \leq  2d(\xi,\eta)
\]
for all $x'$, $x''\in Y$, $x'\neq x''$, $\xi\in \{x',x''\}$, $y\in Y\setminus B(x',2d(x',x''))$, $\eta\in Y\cap B(y,d(x',x'')/8)$. Hence, Lemma \ref{lem:cominur} (iii) implies that
\begin{eqnarray}\label{lem:sxytam6}
\lefteqn{
\int_{Y\cap B(y,d(x',x'')/8)}\frac{d\nu(\eta)}{d(\xi,\eta)^{s_1}d(\eta,y)^{s_2}}
}
\\ \nonumber
&&\qquad
\leq 2^{s_1} d(\xi,y)^{-s_1}\int_{Y\cap B(y,d(x',x'')/8)}\frac{d\nu(\eta)}{d(y,\eta)^{s_2}}
\\ \nonumber
&&\qquad
\leq 2^{s_1} d(\xi,y)^{-s_1}c''_{s_2,Y,Y}8^{-(\upsilon_Y-s_2)}d(x',x'')^{\upsilon_Y-s_2}
\,,
\end{eqnarray}
 for all $x'$, $x''\in Y$, $x'\neq x''$, $\xi\in \{x',x''\}$, $y\in Y\setminus B(x',2d(x',x''))$. Next we note that 
 \[
 A(x',x'',a)\cap\{\eta\in Y: d(\xi,\eta)\geq d(\xi,y)/2\}\subseteq Y\cap B(y,3ad(x',x''))
 \]
 for all $x'$, $x''\in Y$, $x'\neq x''$, $\xi\in \{x',x''\}$, $y\in Y\setminus B(x',2d(x',x''))$. Hence, 	again Lemma \ref{lem:cominur} (iii) implies that
  \begin{eqnarray}\label{lem:sxytam7}
\lefteqn{\int_{A(x',x'',a)\cap\{\eta\in Y: d(\xi,\eta)\geq d(\xi,y)/2\}}\frac{d\nu(\eta)}{d(\xi,\eta)^{s_1}d(\eta,y)^{s_2}}
}
\\ \nonumber
&&\qquad
\leq 2^{s_1}d(\xi,y)^{-s_1}\int_{Y\cap B( y,3ad(x',x''))}\frac{d\nu(\eta)}{d(\eta,y)^{s_2}}
\\ \nonumber
&&\qquad
\leq 2^{s_1}d(\xi,y)^{-s_1}c''_{s_2,Y,Y}(3a)^{\upsilon_Y-s_2}d(x',x'')^{\upsilon_Y-s_2}
\end{eqnarray}
 and that
 \begin{eqnarray}\label{lem:sxytam8}
\lefteqn{
 \int_{A(x',x'',a)\cap\{\eta\in Y: d(y,\eta)\geq d(\xi,y)/2\}}\frac{d\nu(\eta)}{d(\xi,\eta)^{s_1}d(\eta,y)^{s_2}}
 }
\\ \nonumber
&&\qquad
\leq 2^{s_2}d(\xi,y)^{-s_2}\int_{Y\cap B(\xi,ad(x',x'')}\frac{d\nu(\eta)}{d(\xi,\eta)^{s_1} }
\\ \nonumber
&&\qquad
\leq 2^{s_2}d(\xi,y)^{-s_2}c''_{	s_1,		Y,Y}a^{\upsilon_Y-s_1}d(x',x'')^{\upsilon_Y-s_1}
 \end{eqnarray}
  for all $x'$, $x''\in Y$, $x'\neq x''$, $\xi\in \{x',x''\}$, $y\in Y\setminus B(x',2d(x',x''))$. Then by combining inequalities (\ref{lem:sxytam4})--(\ref{lem:sxytam8}) and Lemma \ref{lem:rec}, we deduce the validity of the inequality (\ref{lem:sxytam1}) of the statement.\hfill  $\Box$ 

\vspace{\baselineskip}

Finally, we now show that Theorem \ref{thm:cokewa} implies the validity   the following statement, that extends a classical result of Giraud \cite{Gi34} (see also Miranda \cite[Chap II, \S 11, p.~24]{Mi70}).
\begin{theorem}\label{thm:cokeaba}
 Let  $Y$  be a bounded subset  of a metric space $(M,d)$. Let   $\nu$ be as in (\ref{eq:nu}) with $X=Y$. Let $\upsilon_Y\in]0,+\infty[$.
Let $Y$ be upper $\upsilon_Y$-Ahlfors regular with respect to $Y$. Let $s_1$, $t_1\in[0,\upsilon_Y[$. 
 Let $s_3\in [0,+\infty[$.  	Let $s_2'\in[0,\upsilon_Y[$, $s_2''\in [0,s_3]$, $s_2\equiv s_2'+s_2''$. 
 If either $s_1+t_1=\upsilon_Y$ or $s_2'+t_1=\upsilon_Y$, then we further assume that $Y$ is strongly upper $\upsilon_Y$-Ahlfors regular with respect to $Y$. 
Let $K_{1} \in {\mathcal{K}}_{s_1, s_2,s_3} (Y\times Y)$, 
$K_2\in  {\mathcal{K}}_{t_1, Y\times Y}$.  Let $K_3$ be the composite kernel of $K_1$ and $K_2$ (cf.~(\ref{thm:cokewa0})). Then the following statements hold.
\begin{enumerate}
\item[(i)] If    $s_1+t_1<\upsilon_Y$, 
$s_2'+t_1<\upsilon_Y$, then the kernel  $K_3$ belongs to the class ${\mathcal{K}}_{0, \max\{s_1,t_1\},\min\{s_3-s_2'',\upsilon_Y-s_1, \upsilon_Y-t_1\}} (Y\times Y)$.
 \item[(ii)] If  $s_1+t_1<\upsilon_Y$, 
$s_2'+t_1=\upsilon_Y$, then the kernel  $K_3$ belongs to the class ${\mathcal{K}}_{0, \max\{\epsilon,s_1,t_1\},\min\{s_3-s_2'',\upsilon_Y-s_1, \upsilon_Y-t_1\}} (Y\times Y)$ for all $\epsilon\in]0,+\infty[$.
 \item[(iii)] If  $s_1+t_1<\upsilon_Y$, 
$s_2'+t_1>\upsilon_Y$, then the kernel $K_3$ belongs to the class  ${\mathcal{K}}_{0,
\max\{s_2'+t_1-\upsilon_Y, s_1,t_1\},\min\{s_3-s_2'',\upsilon_Y-s_1, \upsilon_Y-t_1\}} (Y\times Y)$.
 \item[(iv)] If  $s_1+t_1=\upsilon_Y$, 
$s_2'+t_1<\upsilon_Y$, then the kernel $K_3$ belongs to the class ${\mathcal{K}}_{\epsilon, \max\{s_1,t_1\},\min\{s_3-s_2'',\upsilon_Y-s_1, \upsilon_Y-t_1\}} (Y\times Y)$  for all $\epsilon\in]0,+\infty[$.
 \item[(v)] If   $s_1+t_1=\upsilon_Y$, 
$s_2'+t_1=\upsilon_Y$, then the kernel $K_3$ belongs to the class ${\mathcal{K}}_{\epsilon, \max\{\epsilon,s_1,t_1\},\min\{s_3-s_2'',\upsilon_Y-s_1, \upsilon_Y-t_1\}} (Y\times Y)$ for all $\epsilon\in]0,+\infty[$.
 \item[(vi)] If   $s_1+t_1=\upsilon_Y$, 
$s_2'+t_1>\upsilon_Y$, then the kernel  $K_3$ belongs to the class ${\mathcal{K}}_{\epsilon, 
\max\{
s_2'+t_1-\upsilon_Y,s_1,t_1\},\min\{s_3-s_2'',\upsilon_Y-s_1, \upsilon_Y-t_1\}} (Y\times Y)$ for all $\epsilon\in]0,+\infty[$.
 \item[(vii)] If   $s_1+t_1>\upsilon_Y$, 
$s_2'+t_1<\upsilon_Y$, then the kernel  $K_3$ belongs to the class $ {\mathcal{K}}_{s_1+t_1-\upsilon_Y, \max\{s_1,t_1\},\min\{s_3-s_2'',\upsilon_Y-s_1, \upsilon_Y-t_1\}} (Y\times Y)$.
 \item[(viii)] If  $s_1+t_1>\upsilon_Y$, 
$s_2'+t_1=\upsilon_Y$, then the kernel $K_3$ belongs to the class ${\mathcal{K}}_{s_1+t_1-\upsilon_Y, \max\{\epsilon,s_1,t_1\},\min\{s_3-s_2'',\upsilon_Y-s_1, \upsilon_Y-t_1\}} (Y\times Y)$ for all $\epsilon\in]0,+\infty[$.
 \item[(ix)]   If   $s_1+t_1>\upsilon_Y$ and 
$s_2'+t_1>\upsilon_Y$, then the kernel   $K_3$ belongs to the class   ${\mathcal{K}}_{s_1+t_1-\upsilon_Y, 
\max\{s_2'+t_1-\upsilon_Y, s_1,t_1\},\min\{s_3-s_2'',\upsilon_Y-s_1, \upsilon_Y-t_1\}} (Y\times Y)$.
 \end{enumerate}
\end{theorem}
{\bf Proof.}  By Remark \ref{rem:bddimfi} with $X=Y$, we have $\nu(Y)<+\infty$.  In all cases (i)--(ix), the membership of $K_3$ in the class 
${\mathcal{K}}_{s, Y\times Y}$ for the choice of $s$ as in the first index is an immediate consequence of Theorem \ref{thm:cokewa}. 
  Then we note that if $x'$, $x''\in Y$,  $x'\neq x''$, 	 $y\in Y\setminus B(x',2d(x',x''))$, then   the triangular inequality, Proposition \ref{prop:xytam}, Lemma \ref{lem:sxytam}, 
   the elementary inclusion
$B(x',2d(x',x''))\subseteq B(x'',3d(x',x''))$  and the elementary inequality 
  $d(y,x')\geq 2d(x',x'')$ that holds whenever $y\in Y\setminus B(x',2d(x',x''))$   imply that
\begin{eqnarray}\label{thm:cokeaba1}
\lefteqn{
\left|\int_{Y}K_{1}(x',\eta)K_{2}(\eta,y)\,d\nu(\eta)
-
\int_{Y}K_{1}(x'',\eta)K_{2}(\eta,y)\,d\nu(\eta)\right|
}
\\ \nonumber
&&\quad 
=
\left|\int_{Y\setminus B(x',2 d(x',x''))}[K_{1}(x',\eta)-K_{1}(x'',\eta)]K_{2}(\eta,y)\,d\nu(\eta)
\right|
\\ \nonumber
&&\quad \quad
+
\left|\int_{B(x',2 d(x',x''))}[K_{1}(x',\eta)-K_{1}(x'',\eta)]K_{2}(\eta,y)\,d\nu(\eta)
\right|
\\ \nonumber
&&\quad 
\leq
\|K_1\|_{{\mathcal{K}}_{s_1, s_2,s_3} (Y\times Y)}\int_{Y\setminus B(x',2 d(x',x''))}\frac{d(x',x'')^{s_3}}{d(x',\eta)^{s_2'+s_2''}}|K_{2}(\eta,y)|\,d\nu(\eta)
\\ \nonumber
&&\quad \quad
+\|K_1\|_{{\mathcal{K}}_{s_1, s_2,s_3} (Y\times Y)}\|K_2\|_{{\mathcal{K}}_{t_1, Y\times Y}}
\biggl\{
\int_{B(x',2 d(x',x''))} \frac{d\nu(\eta)}{d(x',\eta)^{s_1}d(\eta,y)^{t_1}}
\\ \nonumber
&&\quad\quad
+
\int_{B(x',2 d(x',x''))} \frac{d\nu(\eta)}{d(x'',\eta)^{s_1}d(\eta,y)^{t_1}}
\biggr\}
\\ \nonumber
&&\quad
\leq
\|K_1\|_{{\mathcal{K}}_{s_1, s_2,s_3} (Y\times Y)}
\|K_2\|_{{\mathcal{K}}_{t_1, Y\times Y}}\int_Y\frac{d(x',x'')^{s_3-s_2''}}{2^{s_2''}d(x',\eta)^{s_2'}d(y,\eta)^{t_1}} \,d\nu(\eta)
\\ \nonumber
&&\quad\quad
+\|K_1\|_{{\mathcal{K}}_{s_1, s_2,s_3} (Y\times Y)}\|K_2\|_{{\mathcal{K}}_{t_1, Y\times Y}}
\biggl\{
\int_{B(x',3 d(x',x''))} \frac{d\nu(\eta)}{d(x',\eta)^{s_1}d(\eta,y)^{t_1}}
\\ \nonumber
&&\quad\quad
+
\int_{B(x'',3 d(x',x''))} \frac{d\nu(\eta)}{d(x'',\eta)^{s_1}d(\eta,y)^{t_1}}
\biggr\}
\\ \nonumber
&&\quad
\leq
2^{-s_2''}d(x',x'')^{s_3-s_2''}\|K\|_{{\mathcal{K}}_{s_1, s_2,s_3} (Y\times Y)}\|K_2\|_{{\mathcal{K}}_{t_1, Y\times Y}}c_{Y,Y,Y}^{(s_2',t_1)}I(x',y)
\\ \nonumber
&&\quad\quad
+2\|K_1\|_{{\mathcal{K}}_{s_1, s_2,s_3} (Y\times Y)}\|K_2\|_{{\mathcal{K}}_{t_1, Y\times Y}}c_{s_1,t_1,3,Y}
\\ \nonumber
&&\quad\quad
\times
\left(d(x',y)^{-t_1}d(x',x'')^{\upsilon_Y-s_1}+d(x',y)^{-s_1}d(x',x'')^{\upsilon_Y-t_1}\right)
\end{eqnarray}
where
\[
I(x',y)
\equiv
\left\{
\begin{array}{ll}
 (1+d(x',y)^{\upsilon_Y-(s_2'+t_1)}) & 
  \text{if}\   s_2'+t_1<\upsilon_Y\,, 
 \\
 (1+|\log d(x',y)|) & 
  \text{if}\ 
 s_2'+t_1=\upsilon_Y\,, 
 \\
  d(x',y)^{\upsilon_Y-(s_2'+t_1)}&
    \text{if}\ 
 s_2'+t_1>\upsilon_Y\,. 
\end{array}
\right.
\]
and we note that
\[
(1+|\log d(x',y)|)=(1+|\log d(x',y)|)d(x',y)^\epsilon\frac{1}{d(x',y)^\epsilon}
\]
and that for each $\epsilon\in]0,+\infty[$, the functions
\begin{eqnarray}\label{thm:cokeaba2}
&&(1+d(x',y)^{\upsilon_Y-(s_2'+t_1)})\quad\text{in\ case}\ s_2'+t_1<\upsilon_Y\,,
\\ \nonumber
&&(1+|\log d(x',y)|)d(x',y)^\epsilon\quad\text{in\ case}\ s_2'+t_1=\upsilon_Y\,,
\end{eqnarray}
are bounded in $(x',y)\in (Y\times Y)\setminus  {\mathbb{D}}_{Y\times Y}$ provided that $Y$ is bounded. Here we note that
our application of Lemma \ref{lem:sxytam} has required that $s_1$, $s_2'\in[0,\upsilon_Y[$ and that $Y$ is strongly upper $\upsilon_Y$ regular when $s_2'+t_1=\upsilon_Y$. Hence we conclude that $K_3$ belongs to the class 
${\mathcal{K}}_{a_1,a_2,a_3}(Y\times Y)$ for the parameters $a_1$, $a_2$, $a_3$ indicated in all statements 
  (i)--(ix).\hfill  $\Box$ 

\vspace{\baselineskip}

\section{A regularity theorem for Fredholm integral e\-quations of the second kind on upper Ahlfors regular sets}\label{sec:corefr}

\begin{theorem}\label{thm:cizere}
 Let  $Y$  be a bounded subset  of a metric space  $(M,d)$.  Let $\nu$ be as in (\ref{eq:nu}) with $X=Y$. Let $\upsilon_Y\in]0,+\infty[$. Let $s\in [0,\upsilon_Y[$.     
Let $Y$ be upper $\upsilon_Y$-Ahlfors regular with respect to $Y$.
Let $K\in {\mathcal{K}}_{s,Y\times Y}$.
Let   $A[K,\cdot]$ be the integral operator in $L^{2}_\nu(Y)$ that  is associated to the kernel $K$ (cf.~(\ref{prop:k0a0})). If $g\in C^{0}_b(Y)$ and if $\mu\in L^{2}_\nu(Y)$ satisfies the inequality
\begin{equation}
\label{thm:cizere1}
\mu-A[K,\mu]=g\qquad\nu-{\mathrm{a.e.\ in}}\ Y\,,
\end{equation}
then   $\mu\in C_b^{0}(Y)$. 
\end{theorem}
{\bf Proof.} In order to shorten our notation, we set $A[\cdot]\equiv A[K,\cdot]$. We first prove by induction that
\begin{equation}
\label{thm:cizere2}
A^{r}[\mu]=\mu-\sum_{j=0}^{r-1}A^{j}[g]\,,
\end{equation}
for all $r\in {\mathbb{N}}\setminus\{0\}$. If $r=1$, we have
\[
A[\mu]=\mu-A^{0}[g]=\mu-g
\]
by assumption, and thus equality (\ref{thm:cizere2}) holds.
We now assume that  (\ref{thm:cizere2}) holds for $r\in {\mathbb{N}}\setminus\{0\}$, and we prove it for $r+1$.
By applying $A$ to both hand sides of  equality (\ref{thm:cizere2}) for $r\in {\mathbb{N}}\setminus\{0\}$, we obtain
\[
A^{r+1}[\mu]=A[\mu]-\sum_{j=0}^{r-1}A^{j+1}[g]=
\mu-g-\sum_{j=1}^{r}A^{j}[g]=\mu= \mu-\sum_{j=0}^{r}A^{j}[g]\,,
\]
and thus equality equality (\ref{thm:cizere2}) holds for $r+1$. Then the induction principle implies that equality (\ref{thm:cizere2}) holds for all $r\in {\mathbb{N}}\setminus\{0\}$.
By Remark \ref{rem:bddimfi} with $X=Y$, we have $\nu(Y)<+\infty$. 
Since $s\in[0,\upsilon_Y[$, and $K$ is continuous outside of the diagonal, it is known that the operator $A[\cdot]$ is linear and continuous from $L^{\infty}_\nu(Y)$ to $ C^{0}_b(Y)$ (cf. \textit{e.g.}, \cite[Prop.~4.1]{La22a}, \cite[Prop.~4.3]{La22b}). 
Then our assumption that  $g\in   L^{\infty}_\nu(Y)$,    implies that 
  $A[g]\in C^{0}_b(Y)$ and similarly that $A^{j}[g]\in C^{0}_b(Y)$ for all $j\in {\mathbb{N}}$. Since $Y$ is bounded and the kernel $K$ has a weak singularity and is continuous on $Y^{2}\setminus {\mathbb{D}}_{	 Y\times Y 		}$, Corollary \ref{corol:iterwsa} ensures that 
 there exists $r\in {\mathbb{N}}\setminus\{0\}$ such that $K^{(r)}$ has a continuous and bounded extension to     $Y^{2}$. 
Then the continuity and boundedness of the extension of $K^{(r)}$ to $Y^2$ and classical theorems of continuity for integrals depending on a parameter imply that $A^{r}[\mu]\in C^{0}_b(Y)$. Then we have
\[
\mu=A^{r}[\mu]+\sum_{j=0}^{r-1}A^{j}[g]\in C^{0}_b(Y)\,.
\]
\hfill  $\Box$ 

\vspace{\baselineskip}

\section{A H\"{o}lder regularity theorem for Fredholm integral equations of the second kind on upper Ahlfors regular sets}\label{sec:horefr}
We now turn to show the (generalized) H\"{o}lder regularity of the solutions of a Fredholm integral equation of the second kind when the datum is H\"{o}lder continuous.
\begin{theorem}\label{thm:cizereh}
  Let  $Y$  be a bounded subset  of a metric space  $(M,d)$.   Let $\nu$ be as in (\ref{eq:nu}) with $X=Y$. Let $\upsilon_Y\in]0,+\infty[$.  Let $Y$ be upper $\upsilon_Y$-Ahlfors regular with respect to $Y$. Let $\theta\in]0,1]$, 
\[
s_1\in [0,\upsilon_Y[\cap [\upsilon_Y-1,\upsilon_Y[\,,\quad
s_2\in [0,+\infty[\,,\quad s_3\in ]0,1]\,.
\]

If $s_{2} = \upsilon_Y$, we further require that   $Y$ be strongly upper $\upsilon_Y$-Ahlfors regular with respect to $X$.

If $s_{2} >\upsilon_Y$, we further require that $s_2<\upsilon_Y+ s_3$. 

Let $\varpi$ be the map from $[0,+\infty[$ to itself defined by $\varpi(0)\equiv 0$ and 
\begin{equation}\label{prop:k0aa0}
\varpi(r)\equiv \left\{
\begin{array}{ll}
r^{\min\{\upsilon_Y-s_1,s_{3}\}}  & \text{if}\ s_{2}<\upsilon_Y\,,
\\
\max\{
r^{\upsilon_Y-s_1},\omega_{s_{3}}(r)\} & \text{if}\ s_{2}=\upsilon_Y\,,
 \\
r^{\min\{
\upsilon_Y-s_1, s_{3}+\upsilon_Y- s_{2} 
\} }&\text{if}\ s_{2} >\upsilon_Y\,,
\end{array}
\right.
\ \ \forall r\in]0,+\infty[\,.
\end{equation}
Let $K\in  {\mathcal{K}}_{ s_{1},s_{2},s_{3}   }(Y\times Y)$.
Let   $A[K,\cdot]$ be the integral operator in $L^{2}_\nu(Y)$ that is associated to the kernel $K$ (cf.~(\ref{prop:k0a0}) with $X=Y$). 

If  $g\in C^{0,\theta}_b(Y)$ and if $\mu\in L^{2}_\nu(Y)$ satisfies the equality
\begin{equation}
\label{thm:cizereh1}
\mu-A[K,\mu]=g\qquad\nu-{\mathrm{a.e.\ in}}\ Y\,,
\end{equation}
then   $\mu\in C_b^{0,\max\{r^\theta,\varpi(r)\}}(Y)$. 
\end{theorem}
{\bf Proof.} By Theorem \ref{thm:cizere}, we know that $\mu\in C_b^{0}(Y)$. By Remark \ref{rem:bddimfi} with $X=Y$, we have $\nu(Y)<+\infty$.  Then under our assumptions, \cite[Prop.~5.2] {La22a} ensures that $A[K,\cdot]$ is linear and continuous from $C_b^{0}(Y)$ to $C_b^{0,\varpi(\cdot)}(Y)$. Hence,
  $A[K,\mu]$ belongs to $C_b^{0,\varpi(\cdot)}(Y)$. Since $g\in C^{0,\theta}_b(Y)$, equality (\ref{thm:cizereh1}) implies that $\mu\in C_b^{0,\max\{r^\theta,\varpi(r)\}}(Y)$.\hfill  $\Box$ 

\vspace{\baselineskip}

\begin{example}
Let $n\in {\mathbb{N}}\setminus\{0,1\}$, $\alpha\in]0,1]$. Let $\Omega$ be a bounded open Lipschitz subset of ${\mathbb{R}}^n$.
If $\alpha=1$ we further require that $\Omega$ is of class $C^1$. Let $m_{n-1}$ be the $(n-1)$-dimensional measure on $\partial\Omega$.   Let 
\[
K\in {\mathcal{K}}_{(n-1)-\alpha,n-\alpha,1}((\partial\Omega)\times(\partial\Omega))\,.
\]
 Let $A[K,\cdot]$ be the integral operator in $L^{2}_\nu(\partial\Omega)$ that is associated to the kernel $K$ (cf.~(\ref{prop:k0a0}) with $X=Y=\partial\Omega$). If  $g\in C^{0,\alpha}_b(\partial\Omega)$   and   $\mu\in L^{2}_\nu(\partial\Omega)$ satisfy the equality
\begin{equation}
\label{thm:cizerehe1}
\mu-A[K,\mu]=g\qquad m_{n-1}-{\mathrm{a.e.\ in}}\ \partial\Omega\,,
\end{equation}
then   $\mu\in C_b^{0,\alpha}(\partial\Omega)$ for $\alpha\in ]0,1[$ and $\mu\in C_b^{0,\omega_1(\cdot)}(\partial\Omega)$ for $\alpha=1$ (cf.~(\ref{omth})).
\end{example}
{\bf Proof.} We first observe that    $Y\equiv \partial\Omega$ is upper $m_{n-1}$-Ahlfors regular with respect to $X\equiv \partial\Omega$
and that  if $\Omega$ is also of class $C^1$, then $\partial\Omega$ is also
strongly upper $m_{n-1}$-Ahlfors regular with respect to $\partial\Omega$ (cf.~\textit{e.g.}, \cite[Prop.~6.5, Rem.~6.2]{La24}) and we can set
\[
\upsilon_Y\equiv n-1\,,\quad s_1\equiv (n-1)-\alpha\,,\quad s_2\equiv n-\alpha\,, \quad s_3\equiv 1\,.
\]
In particular, $s_2=n-\alpha<\upsilon_Y+s_3=(n-1)+1$ if $n-\alpha=s_2>(n-1)=\upsilon_Y$, \textit{i.e.}, if $\alpha\in]0,1[$.
Then Theorem \ref{thm:cizereh} ensures that $\mu\in C_b^{0,\max\{r^\alpha,\varpi(r)\}}(\partial\Omega)$. If $\alpha\in]0,1[$, then
$s_2\equiv n-\alpha=(n-1)+(1-\alpha)>n-1$ and
\[
C_b^{\varpi(r)}(\partial\Omega)=C_b^{r^{\min\{
(n-1)-[(n-1)-\alpha], 1+(n-1)- [(n-1)+(1-\alpha)] 
\} }}(\partial\Omega)=C_b^{\alpha}(\partial\Omega)\,.
\]
If instead $\alpha=1$, then
$
s_2\equiv n-\alpha=(n-1)+(1-1)=n-1 
$
and
\[
C_b^{\varpi(r)}(\partial\Omega)=C_b^{
\max\{
r^{(n-1)-[(n-1)-1]},\omega_{1}(r)\}
}(\partial\Omega)=C_b^{\omega_{1}(r)}(\partial\Omega) 
\]
and thus the proof is complete.\hfill  $\Box$ 

\vspace{\baselineskip}

We note that under the assumptions of the previous example, if $\Omega$ is of class $C^{1,\alpha}$, then the kernel of the double layer potential that is associated to a fundamental solution of a second order elliptic operator with constant coefficients is known to be of class ${\mathcal{K}}_{(n-1)-\alpha,n-\alpha,1}((\partial\Omega)\times(\partial\Omega))$ (see reference \cite[Rem. 6.1]{DoLa17} with Dondi).

We now show that if we have some information on the H\"{o}lder regularity of  $A[K,1]$, then we can prove more information on the (generalized) H\"{o}lder continuity of the solution
of a Fredholm integral equation of the second kind when the datum is H\"{o}lder continuous.
\begin{theorem}\label{thm:cizereht}
Let  $Y$  be a bounded subset  of a metric space  $(M,d)$.   Let $\nu$ be as in (\ref{eq:nu}) with $X=Y$. Let $\upsilon_Y\in]0,+\infty[$. 
   Let $\theta\in]0,1]$. Let
\[
s_1\in [0,\upsilon_Y[\,,\quad
  \beta\in]0,1]\,,\quad 
 s_{2}\in [\beta,+\infty[\,,\quad
 s_{3}\in]0,1]\,.
\]
Let $Y$ be upper $\upsilon_Y$-Ahlfors regular with respect to $Y$.
If   $s_{2}-\beta= \upsilon_Y$, we further require that   $Y$ be strongly upper $\upsilon_Y$-Ahlfors regular with respect to $Y$.

If $s_{2}-\beta>\upsilon_Y$, we further require that $s_{3}+\upsilon_Y-(s_{2}-\beta)>0$.

 Let $\omega$ be the map from $[0,+\infty[$ to itself defined by $\omega(0)\equiv 0$  
 and 
\begin{equation}\label{prop:k0b1aa}
\omega(r)\equiv \left\{
\begin{array}{ll}
r^{\min\{\upsilon_Y-s_1+\beta,s_{3}\}}  & \text{if}\ s_{2}-\beta<\upsilon_Y\,,
\\
\max\{
r^{\upsilon_Y-s_1+\beta},\omega_{s_{3}}(r)\}  & \text{if}\ s_{2}-\beta=\upsilon_Y\,,
 \\
r^{\min\{
\upsilon_Y-s_1+\beta, s_{3}+\upsilon_Y-(s_{2}-\beta) 
\} }&\text{if}\ s_{2}-\beta>\upsilon_Y\,,
\end{array}
\right.
\qquad\forall r\in]0,+\infty[\,.
\end{equation}
Let $K\in{\mathcal{K}}_{ s_{1},s_{2},s_{3}   }(Y\times Y)$.
Let   $A[K,\cdot]$ be the integral operator in $L^{2}_\nu(Y)$ that is associated to the kernel $K$ (cf.~(\ref{prop:k0a0}) with $X=Y$). Let 
\begin{equation}\label{thm:cizereht2}
A[K,1]\in C^{0,\omega(\cdot)}_b(Y)\,.
\end{equation}
If $g\in C^{0,\theta}_b(Y)$ and if $\mu\in C^{0,\beta}_b (Y)$ satisfies the equality
\begin{equation}
\label{thm:cizereht3}
\mu-A[K,\mu]=g\qquad\nu-{\mathrm{a.e.\ in}}\ Y\,.
\end{equation}
Then   $\mu\in C_b^{0,\max\{r^\theta,\omega(r)\}}(Y)$. 
\end{theorem}
{\bf Proof.} By assumption, we know that $\mu\in C_b^{0,\beta}(Y)$.  By assumption (\ref{thm:cizereht2}), we have 
$A[K,1]\in C^{0,\omega(\cdot)}_b(Y)$. By Remark \ref{rem:bddimfi}  with $X=Y$, we have $\nu(Y)<+\infty$. Then \cite[Prop.~5.11]{La22a}
  and Remark \ref{rem:om4} imply that $A[K,\mu]\in C_b^{0,\omega(\cdot)}(Y)$. Since $g\in C^{0,\theta}_b(Y)$, equality (\ref{thm:cizereht3}) implies that $\mu\in C_b^{0,\max\{r^\theta,\omega(r)\}}(Y)$.\hfill  $\Box$ 

\vspace{\baselineskip}

 We note that under the assumptions of Theorem \ref{thm:cizereht}, if $C_b^{0,\omega(\cdot)}(Y)\subseteq C_b^{0,\theta}(Y)$, then we can infer that $\mu$ is actually of class $C_b^{0,\theta}(Y)$ as the datum $g$. 

\vspace{\baselineskip}

 \noindent
{\bf Statements and Declarations}\\

 \noindent
{\bf Competing interests:} This paper does not have any  conflict of interest or competing interest.

 \noindent
{\bf Acknowledgement.}  
This paper represents an extension of   the work performed 
by M.~Norman in his `Laurea Triennale' dissertation \cite{No25} under the guidance of 
M.~Lanza de Cristoforis. 
The first author  acknowledges  the support of GNAMPA-INdAM,   of the Project funded by the European Union – Next Generation EU under the National Recovery and Resilience Plan (NRRP), Mission 4 Component 2 Investment 1.1 - Call for tender PRIN 2022 No. 104 of February, 2 2022 of Italian Ministry of University and Research; Project code: 2022SENJZ3 (subject area: PE - Physical Sciences and Engineering) ``Perturbation problems and asymptotics for elliptic differential equations: variational and potential theoretic methods''  and the support from the EU through the H2020-MSCA-RISE-2020 project EffectFact, Grant agreement ID: 101008140.

.

\end{document}